\documentclass[a4paper,10pt]{article}
\usepackage[affil-it]{authblk}

\usepackage[english]{babel}

\usepackage{mathrsfs}
\usepackage{amssymb}
\usepackage{amsmath}
\usepackage{latexsym}
\usepackage{amsthm}
\usepackage{amscd}
\usepackage{graphicx}

\usepackage[colorlinks=true]{hyperref}

\theoremstyle{plain}
\newtheorem{teo}{Theorem}[section]
\newtheorem{prop}[teo]{Proposition}
\newtheorem{coroll}[teo]{Corollary}
\newtheorem{lem}[teo]{Lemma}

\theoremstyle{definition}
\newtheorem{defin}[teo]{Definition}
\newtheorem{ese}[teo]{Example}
\newtheorem{rmk}[teo]{Remark}

\numberwithin{equation}{section}

\newcommand{\R}{\mathbb R}
\newcommand{\Ha}{\mathcal H}
\newcommand{\Ll}{\mathcal L}
\newcommand{\Co}{\mathcal C}
\newcommand{\eps}{\varepsilon}
\newcommand{\Dim}{\textrm{Dim}}

\newcommand{\s}{\mathbb S}

\title{Fractional perimeter from a fractal perspective}

\author{Luca Lombardini}

\affil{\footnotesize Universit\`a degli Studi di Milano\\ Via Cesare Saldini, 50\\ 20133, Milano, Italia}

\date{}

\begin{document}
\maketitle

\begin{abstract}
%It is well known that open sets with a regular boundary have finite fractional perimeter.
%In this paper we show that sets with an irregular, ``fractal'', boundary can have finite fractional perimeter.\\
Following \cite{Visintin}, we exploit the fractional perimeter of a set to give a definition of fractal dimension
for its measure theoretic boundary.\let\thefootnote\relax\footnotetext{I would like to express my gratitude to Prof. Enrico Valdinoci
for his advice, support and patience.}

We calculate the fractal dimension of sets which can be defined in a recursive way and
we give some examples of this kind of sets, explaining how to construct them starting from well known self-similar fractals.\\
In particular, we show that in the case of the von Koch snowflake $S\subset\R^2$ this fractal dimension
coincides with the Minkowski dimension, namely
\begin{equation*}
P_s(S)<\infty\qquad\Longleftrightarrow\qquad s\in\Big(0,2-\frac{\log4}{\log3}\Big).
\end{equation*}

We also study the asymptotics as $s\to1^-$ of the fractional perimeter %$P_s(E,\Omega)$
of a set having finite (classical) perimeter.
%in a bounded open set with Lipshchitz boundary.

\end{abstract}

\tableofcontents

\begin{section}{Introduction and main results}

It is well known (see e.g. \cite{Gamma} and \cite{cafenr}) that sets with a regular boundary have finite $s$-fractional perimeter
for every $s\in(0,1)$.

In this paper we show that also sets with an irregular, ``fractal'', boundary can have finite $s$-perimeter for every
$s$ below some threshold $\sigma<1$.\\
Actually, the $s$-perimeter can be used to define a ``fractal dimension'' for the measure theoretic boundary
\begin{equation*}
\partial^-E:=\{x\in\R^n\,|\,0<|E\cap B_r(x)|<\omega_nr^n\textrm{ for every }r>0\},
\end{equation*}
of a set $E\subset\R^n$.
Indeed, in \cite{Visintin} the author suggested using the index $s$ of the seminorm $[\chi_E]_{W^{s,1}}$ as a way to measure the codimension
of $\partial^-E$
%, defining
%\begin{equation*}
%\Dim_F(\partial^-E,\Omega):=n-\sup\{s\in(0,1)\,|\,[\chi_E]_{W^{s,1}(\Omega)}<\infty\}.
%\end{equation*}
and he proved that the fractal dimension obtained in this way is less or equal than the (upper) Minkowski dimension.

We give an example of a set, the von Koch snowflake, for which these two dimensions coincide.

%This represents a deep difference between the fractional and the classical perimeter. Having (locally) finite $s$-perimeter
%is not a condition strong enough to guarantee an analogue of De Giorgi's structure Theorem

Moreover, exploiting the roto-translation invariance and the scaling property of the $s$-perimeter,
we calculate the dimension
of sets which can be defined in a recursive way similar to that of the von Koch snowflake.

%Since the $s$-perimeters of two sets which differ only in a set of measure zero are the same,
%this fractal dimension is not suited to detect the fractal nature of sets of measure zero, like the Sierpinski triangle.\\
%Indeed, the fractal object

%Still, we show how to explit the structure of the Sierpinski triangle and other self siilar fractals, to define

On the other hand, as remarked above, sets with a regular boundary have finite $s$-perimeter for every $s$ and actually
their $s$-perimeter converges, as $s$ tends to 1, to the classical perimeter,
both in the classical sense (see \cite{cafenr}) and in the $\Gamma$-convergence sense (see \cite{Gamma}).\\
As a simple byproduct of the computations developed in this paper, 
we exploit Theorem 1 of \cite{Davila} to prove this asymptotic property for a set $E$ having finite classical
perimeter in a bounded open set with Lipschitz boundary.\\
This last result is probably well known to the expert, though not explicitly stated in the literature (as far as we know).\\
In particular, we remark that this lowers the regularity
requested in \cite{cafenr}, where the authors asked the boundary $\partial E$ to be $C^{1,\alpha}$.\\

We begin by recalling the definition of $s$-perimeter.

Let $s\in(0,1)$ and let $\Omega\subset\mathbb R^n$ be an open set. The $s$-fractional perimeter
of a set $E\subset\mathbb R^n$ in $\Omega$ is defined as
\begin{equation*}
P_s(E,\Omega):=\mathcal L_s(E\cap\Omega,\Co E\cap\Omega)+
\mathcal L_s(E\cap\Omega,\Co E\setminus\Omega)+
\mathcal L_s(E\setminus\Omega,\Co E\cap\Omega),
\end{equation*}
where
\begin{equation*}
\mathcal L_s(A,B):=\int_A\int_B\frac{1}{|x-y|^{n+s}}\,dx\,dy,
%=\int_{\R^n}\int_{\R^n}\frac{1}{|x-y|^{n+s}}\chi_A(x)\chi_B(y)\,dx\,dy
\end{equation*}
for every couple of disjoint sets $A,\,B\subset\mathbb R^n$. We simply write $P_s(E)$ for $P_s(E,\R^n)$.\\
%Notice that formally this coincides with
%\begin{equation}\label{formally_perimeter}
%P_s(E,\Omega)=\frac{1}{2}\big([\chi_E]_{H^\frac{s}{2}(\mathbb R^n)}^2-[\chi_E]_{H^\frac{s}{2}(\Omega^c)}^2\big),
%\end{equation}
%where $[u]_{H^\sigma(\mathcal O)}$ denotes the Gagliardo seminorm
%\begin{equation*}
%[u]_{H^\sigma(\mathcal O)}=\Big(\int_{\mathcal O}\int_{\mathcal O}\frac{|u(x)-u(y)|^2}{|x-y|^{n+2\sigma}}\,dx\,dy\Big)^\frac{1}{2}.
%\end{equation*}
We can also write the fractional perimeter as the sum
\begin{equation*}
P_s(E,\Omega)=P_s^L(E,\Omega)+P_s^{NL}(E,\Omega),
\end{equation*}
where
\begin{equation*}\begin{split}
&P_s^L(E,\Omega):=\mathcal L_s(E\cap\Omega,\Co E\cap\Omega)=\frac{1}{2}[\chi_E]_{W^{s,1}(\Omega)},\\
&
P_s^{NL}(E,\Omega):=\Ll_s(E\cap\Omega,\Co E\setminus\Omega)+\Ll_s(E\setminus\Omega,\Co E\cap\Omega).
\end{split}\end{equation*}
We can think of $P^L_s(E,\Omega)$ as the local part of the fractional perimeter, in the sense that if $|(E\Delta F)\cap\Omega|=0$,
then $P^L_s(F,\Omega)=P^L_s(E,\Omega)$.

We say that a set $E$ has locally finite $s$-perimeter if it has finite $s$-perimeter in
every bounded open set $\Omega\subset\R^n$.\\

Now we give precise statements of the results obtained,
starting with the fractional analysis of fractal dimensions.

\begin{subsection}{Fractal boundaries}

%A natural question is: what kind of sets have finite fractional perimeter?

%If the set is regular enough, then it has finite $s$-perimeter for every $s\in(0,1)$ (see Section 1).

%On the other hand, we show that also sets having a fractal boundary can have finite $s$-perimeter for $s$ below some threshold
%$\sigma\in(0,1)$.

%Indeed, following
%\cite{Visintin}, we exploit the $s$-perimeter to give a notion of fractal dimension for the measure theoretic boundary
%\begin{equation*}
%\partial^-E:=\{x\in\R^n\,|\,0<|E\cap B_r(x)|<\omega_nr^n\textrm{ for every }r>0\}.
%\end{equation*}
%(see Remark $\ref{gmt_assumption}$ below or Appendix C for the formal definition).

First of all, we prove in Section 3.1 that in some sense the measure theoretic boundary $\partial^-E$ is the ``right definition'' of boundary for
working with the $s$-perimeter.

To be more precise, we show that
\begin{equation*}
\partial^-E=\{x\in\R^n\,|\,P_s^L(E,B_r(x))>0,\,\forall\,r>0\},
\end{equation*}
and that if $\Omega$ is a connected open set, then
\begin{equation*}
P_s^L(E,\Omega)>0\quad\Longleftrightarrow\quad \partial^-E\cap\Omega\not=\emptyset.
\end{equation*}
This can be thought of as an analogue in the fractional framework of the fact that for a Caccioppoli set $E$ we have $\partial^-E=$ supp $|D\chi_E|$.
%We show that something similar holds also in this fractional framework.\\

Now the idea of the definition of the fractal dimension consists in using the index $s$ of $P_s^L(E,\Omega)$ to measure the codimension of
$\partial^- E\cap\Omega$,
\begin{equation*}
\Dim_F(\partial^-E,\Omega):=n-\sup\{s\in(0,1)\,|\,P^L_s(E,\Omega)<\infty\}.
\end{equation*}

As shown in \cite{Visintin} (Proposition 11 and Proposition 13), the fractal dimension $\textrm{Dim}_F$ defined in this way is related to the (upper) Minkowski dimension by
\begin{equation}\label{intro_dim_ineq}
\Dim_F(\partial^-E,\Omega)\leq\overline{\Dim}_\mathcal M(\partial^-E,\Omega),
\end{equation}
(for the convenience of the reader we provide a proof in 
Proposition $\ref{vis_prop}$).

If $\Omega$ is a bounded open set with Lipschitz boundary, % (or $\Omega=\R^n$ if $E$ is bounded), 
this means that
\begin{equation}\label{intro_dim_ineq2}
P_s(E,\Omega)<\infty\qquad\textrm{for every }s\in\big(0,n-\overline{\Dim}_\mathcal M(\partial^-E,\Omega)\big),
\end{equation}
%when $\Omega$ is a bounded open set with Lipschitz boundary.
%We remark that if $\Omega$ is bounded and has Lipschitz boundary, then 
since the nonlocal part of the $s$-perimeter of any set $E\subset\R^n$ is
\begin{equation*}
P_s^{NL}(E,\Omega)\leq2P_s(\Omega)<\infty,\qquad\textrm{for every }s\in(0,1).
\end{equation*}

%  (see Corollary
%$\ref{embedding_fin_per_coroll}$)
%so in this case using just.\\

%However in general we do not know what happens for $s$ above the threshold appearing in $(\ref{intro_dim_ineq2})$.

We show that for the von Koch snowflake $(\ref{intro_dim_ineq})$
is actually an equality.%this is sharp.

\begin{figure}[htbp]
\begin{center}
\includegraphics[width=100mm]{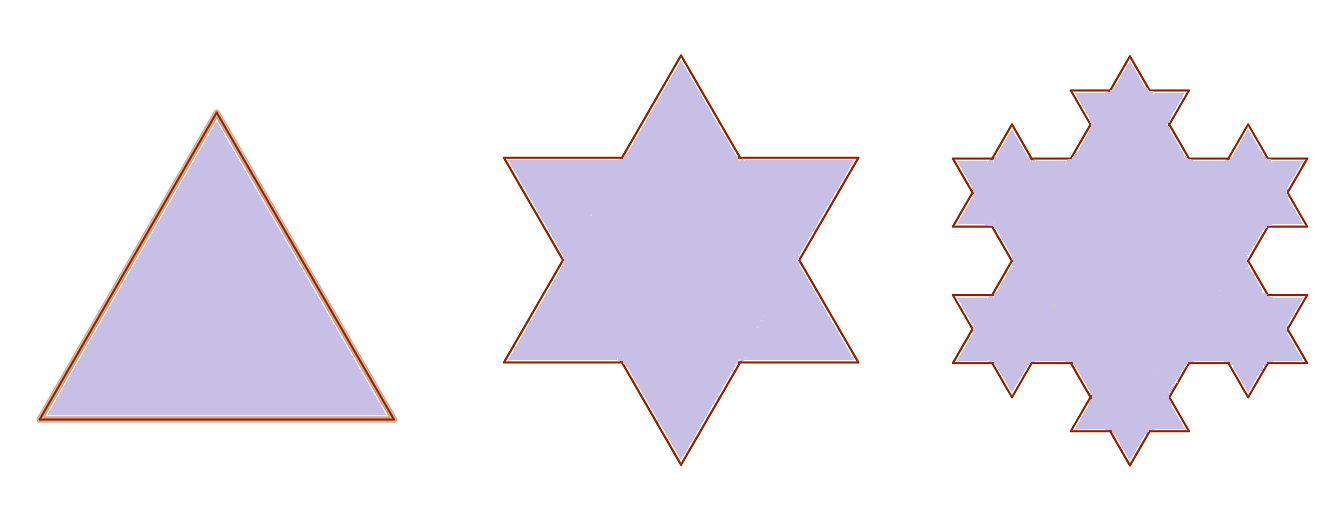}
\caption{{\it The first three steps of the construction of the von Koch snowflake}}
\end{center}
\end{figure}

Namely, we prove the following

\begin{teo}[Fractal dimension of the von Koch snowflake]\label{von_koch_snow}
Let $S\subset\R^2$ be the von Koch snowflake. Then
\begin{equation}\label{koch1}
P_s(S)<\infty,\qquad\forall\,s\in\Big(0,2-\frac{\log4}{\log3}\Big),
\end{equation}
and
\begin{equation}\label{koch2}
P_s(S)=\infty,\qquad\forall\,s\in\Big[2-\frac{\log4}{\log3},1\Big).
\end{equation}
Therefore
\begin{equation*}
\Dim_F(\partial S)=\Dim_\mathcal{M}(\partial S)=\frac{\log4}{\log3}.
\end{equation*}
\end{teo}

Actually, exploiting the self-similarity of the von Koch curve, we have
\begin{equation*}
\Dim_F(\partial S,\Omega)=\frac{\log4}{\log3},
\end{equation*}
for every $\Omega$ s.t. $\partial S\cap\Omega\not=\emptyset$.
In particular, this is true for every $\Omega=B_r(p)$ with $p\in S$ and $r>0$ as small as we want.\\

We remark that this represents a deep difference between the classical and the fractional perimeter.\\
Indeed, if a set $E$ has (locally) finite perimeter, then by De Giorgi's structure Theorem we know that its reduced boundary $\partial^*E$
is locally $(n-1)$-rectifiable. Moreover $\overline{\partial^*E}=\partial^-E$, so the reduced boundary is, in some sense,
a ``big'' portion of the measure theoretic boundary.

On the other hand, there are (open) sets, like the von Koch snowflake, which have a ``nowhere rectifiable'' boundary
(meaning that $\partial^-E\cap B_r(p)$ is not $(n-1)$-rectifiable for every $p\in\partial^-E$ and $r>0$)
and still have finite $s$-perimeter for every $s\in(0,\sigma_0)$.\\

Moreover our argument for the von Koch snowflake is quite general and can be adapted to
%show that $(\ref{snowflake_intro_dim})$ actually holds for
calculate the dimension $\Dim_F$ of
all sets which can be constructed in a similar recursive way 
%similar to that of the snowflake 
(see Section 3.4).\\
Roughly speaking, these sets are defined by adding scaled copies of a fixed ``building block'' $T_0$,
that is
\begin{equation*}
T:=\bigcup_{k=1}^\infty \bigcup_{i=1}^{ab^{k-1}}T_k^i,
\end{equation*}
where $T_k^i:=F_k^i(T_0)$ is a roto-traslation of the scaled set $\lambda^{-k}T_0$ (see Figure 2 below for an example).
We also assume
that $\frac{\log b}{\log\lambda}\in(n-1,n)$.

Theorem $\ref{fractal_bdary_selfsim_dim}$ shows that 
if such a set $T$ satisfies an additional assumption,
namely that ``near'' each set $T_k^i$ we can find a set $S_k^i=F^i_k(S_0)$ contained in $\Co T$,
then the fractal dimension of its measure theoretic boundary is
\begin{equation}
\Dim_F(\partial^-T)=\frac{\log b}{\log\lambda}.
\end{equation}

%Namely this additional assumption is needed to guarantee that the measure theoretic boundary $\partial^-T$
%is not empty and actually .\\

Many well known self-simlar fractals can be written either as (the boundary of) a set $T$ defined as above,
like the von Koch snowflake, or as the difference $E=T_0\setminus T$,
like the Sierpinski triangle and the Menger sponge.

However sets of this second kind are often s.t. $|T\Delta T_0|=0$.\\
Since the $s$-perimeters
of two sets which differ only in a set of measure zero are equal, 
%As a consequence 
in this case the $s$-perimeter can not detect the ``fractal nature'' of $T$.

%of sets of measure zero, like the Sierpinski triangle or the Menger sponge.

Consider for example the Sierpinski triangle, which is defined as $E=T_0\setminus T$ with $T_0$ an equilateral triangle.\\
Then $\partial^-T=\partial T_0$ and
$P_s(T,\Omega)=P_s(T_0,\Omega)<\infty$ for every $s\in(0,1)$.
%and indeed .

%notices only the fractal nature of the measure theoretic boundary, not of the
%topological boundary.

%We remark that many classical examples of self-similar fractal sets, like the Sierpinski triangle and the Menger sponge,
%, do not satisfy this condition.
%In these cases,
%the problem is that
%, so $P_s^L(E,\Omega)=0$. Still, the fractal nature of $E$ is ``contained'' also in $\partial T$, so we can think
%of measuring $P_s(T,\Omega)$ instead.\\
%But this doesn't take us any further.

%Indeed the set $T_0$ is usually a bounded open set with Lipschitz boundary and hence, since
%$|T\Delta T_0|=0$, we get $P_s(T,\Omega)=P_s(T_0,\Omega)$, which is finite for every $s\in(0,1)$.\\

Roughly speaking, the reason of this situation is that the fractal object is the topological boundary of $T$,
while its measure theoretic boundary is regular and has finite (classical) perimeter.\\

%For example, in the case of the Sierpinski triangle
%we have $\partial^-T=\partial T_0$, an equilateral triangle, even if $\partial T$ is a fractal.\\
%As a consequence, the $s$-perimeter does not notice the fractal nature of the set.\\

Still, we show how to modify such self-similar sets, without altering their ``structure'', to obtain new sets
which satisfy the hypothesis of Theorem $\ref{fractal_bdary_selfsim_dim}$. However, the measure theoretic boundary of such
a new set will look quite different from the original fractal (topological) boundary and in general
it will be a mix of smooth parts and unrectifiable parts.
%``self-similar fractal boundaries''.
%We show how to do so for the Sierpinski triangle.

\begin{figure}[htbp]
\begin{center}
\includegraphics[width=60mm]{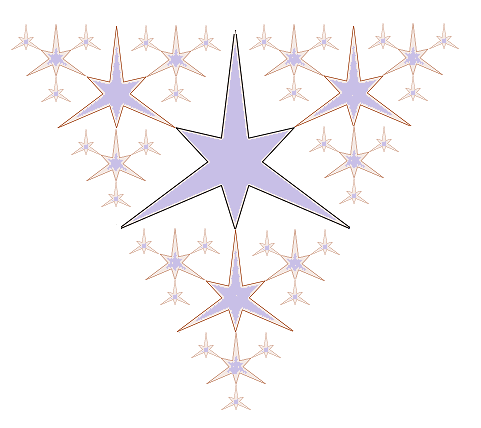}
\caption{{\it Example of a ``fractal'' set constructed exploiting the structure of the Sierpinski triangle (seen at the fourth iterative step), which satisfies the hypothesis of Theorem $\ref{fractal_bdary_selfsim_dim}$}}
\end{center}
\end{figure}

The most interesting examples of 
%sets satisfying the hypothesis of Theorem $\ref{fractal_bdary_selfsim_dim}$
this kind of sets
are probably represented by bounded sets, like the one in Figure 2, because in this case the measure theoretic boundary
does indeed have, in some sense, a ``fractal nature''.\\% (see Remark $\ref{self_sim_frac_bdry_nat_rmk}$).\\
Indeed, if $T$ is bounded, then its boundary $\partial^-T$ is compact. Nevertheless, it has infinite (classical) perimeter
and actually $\partial^-T$ has Minkowski dimension strictly greater than $n-1$, thanks to $(\ref{intro_dim_ineq})$.\\

However, even unbounded sets can have an interesting behavior. Indeed we obtain the following
\begin{prop}\label{expl_farc_prop1}
Let $n\geq2$. For every $\sigma\in(0,1)$ there exists a Caccioppoli set $E\subset\R^n$%, , which has locally finite perimeter
%(and hence also locally finite $s$-perimeter for every $s\in(0,1)$)
%and
s.t.
\begin{equation*}
P_s(E)<\infty\qquad\forall\,s\in(0,\sigma)\quad\textrm{and}\quad P_s(E)=\infty\qquad\forall\,s\in[\sigma,1).
\end{equation*}
\end{prop}

%\begin{rmk}

\noindent
Roughly speaking, the interesting thing about this Proposition is the following. Since $E$ has locally finite perimeter, $\chi_E\in BV_{loc}(\R^n)$,
it also has locally finite $s$-perimeter for every $s\in(0,1)$,
but 
%On the other hand,
the global perimeter $P_s(E)$ is finite if and only if $s<\sigma<1$.
%\end{rmk}

\end{subsection}

\begin{subsection}{Asymptotics as $s\to1^-$}

We have shown that sets with an irregular, eventually fractal, boundary can have finite $s$-perimeter.

On the other hand, if the set $E$ is ``regular'', then it has finite $s$-perimeter for every $s\in(0,1)$.\\
Indeed, if $\Omega\subset\R^n$ is a bounded open set with Lipschitz boundary (or $\Omega=\R^n$),
then $BV(\Omega)\hookrightarrow W^{s,1}(\Omega)$. As a consequence of this embedding,
we obtain
\begin{equation}
P(E,\Omega)<\infty\qquad\Longrightarrow\qquad P_s(E,\Omega)<\infty\quad\textrm{for every }s\in(0,1).
\end{equation}

Actually we can be more precise and obtain a sort of converse,
using only the local part of the $s$-perimeter and adding the condition
\begin{equation*}
\liminf_{s\to1^-}(1-s)P^L_s(E,\Omega)<\infty.
\end{equation*}

Indeed one has the following result, which is just a combination of Theorem
 3' of \cite{BBM} and Theorem 1 of \cite{Davila}, restricted to characteristic functions,
\begin{teo}\label{Davila_conv_local}
Let $\Omega\subset\R^n$ be a bounded open set with Lipschitz boundary. Then
$E\subset\R^n$ has finite perimeter in $\Omega$ if and only if $P_s^L(E,\Omega)<\infty$
for every $s\in(0,1)$, and
\begin{equation}\label{asymptotics_fin_cond}
\liminf_{s\to1}(1-s)P_s^L(E,\Omega)<\infty.
\end{equation}
In this case we have
\begin{equation}\label{asymptotics_local_part}
\lim_{s\to1}(1-s)P_s^L(E,\Omega)=\frac{n\omega_n}{2}K_{1,n}P(E,\Omega).
\end{equation}
\end{teo}
\noindent
We briefly show how to get this result (and in particular why the constant looks like that) from
the two Theorems cited above.

We compute the constant $K_{1,n}$ in an elementary way, showing that
\begin{equation}
\frac{n\omega_n}{2}K_{1,n}=\omega_{n-1}.
\end{equation}

Moreover we show the following
\begin{rmk}
Condition $(\ref{asymptotics_fin_cond})$ is necessary. Indeed, there exist bounded sets (see the following Example) having finite $s$-perimeter for every $s\in(0,1)$ which do not have finite perimeter.\\
This also shows that in general the inclusion $BV(\Omega)\subset\bigcap_{s\in(0,1)}W^{s,1}(\Omega)$
is strict.
\end{rmk}

\begin{ese}\label{inclusion_counterexample}
Let $0<a<1$ and consider the open intervals $I_k:=(a^{k+1},a^k)$ for every $k\in\mathbb{N}$.
Define $E:=\bigcup_{k\in\mathbb{N}}I_{2k}$, which is a bounded (open) set.\\
Due to the infinite number of jumps $\chi_E\not\in BV(\mathbb{R})$. However it can be proved that
$E$ has finite $s$-perimeter for every $s\in(0,1)$. We postpone the proof to Appendix A.
\end{ese}

The main result of Section 2 is the following Theorem, which extends the asymptotic convergence
of $(\ref{asymptotics_local_part})$ to the whole $s$-perimeter, at least when
the boundary $\partial E$ intersects
the boundary of $\Omega$ ``transversally''.

\begin{teo}[Asymptotics]\label{asymptotics_teo}
Let $\Omega\subset\R^n$ be a bounded open set with Lipschitz boundary.
Suppose that $E$ has finite perimeter in $\Omega_\beta$, for some $\beta\in(0,r_0)$, with $r_0>0$ small enough. Then
\begin{equation}\label{asymptotics_nonlocal_estimate}
\limsup_{s\to1}(1-s)P_s^{NL}(E,\Omega)
\leq2\omega_{n-1}\lim_{\rho\to0^+}P(E,N_\rho(\partial\Omega)).
\end{equation}
In particular, if $P(E,\partial\Omega)=0$, then
\begin{equation}
\lim_{s\to1}(1-s)P_s(E,\Omega)=\omega_{n-1}P(E,\Omega).
\end{equation}
Moreover, there exists a set $S\subset(-r_0,\beta)$, at most countable,
s.t.
\begin{equation}\label{asymptotics_ae_convergence}
\lim_{s\to1}(1-s)P_s(E,\Omega_\delta)=\omega_{n-1}P(E,\Omega_\delta),
\end{equation}
for every $\delta\in(-r_0,\beta)\setminus S$.

\end{teo}

Roughly speaking, the second part of this Theorem says that even if we do not have the asymptotic convergence of the $s$-perimeter in $\Omega$, we can slighltly enlarge or restrict $\Omega$ to obtain it.
Actually, since $S$ has null measure, we can restrict or enlarge $\Omega$ as little as we want.\\

In \cite{cafenr} the authors obtained a similar result
for $\Omega=B_R$ a ball, but asking $C^{1,\alpha}$ regularity of $\partial E$
in $B_R$. They proved the convergence in every ball $B_r$ with $r\in(0,R)\setminus S$,
with $S$ at most countable,
exploiting uniform estimates.\\
On the other hand, asking $E$ to have finite perimeter in a neighborhood (as small as we want) of the open set $\Omega$
is optimal.\\
%Indeed if $E\subset\R^n$ is s.t. $(\ref{asymptotics_ae_convergence})$ holds true, then Theorem $\ref{Davila_conv_local}$ guarantees that $E$ has finite perimeter
%in $\Omega_\delta$.\\

In \cite{Gamma} the authors studied the asymptotics as $s\longrightarrow1^-$ in the $\Gamma$-convergence sense.
In particular, for the proof of a $\Gamma$-limsup inequality, which is typically constructive and by density, they show that if $\Pi$ is a polyhedron, then
\begin{equation*}
\limsup_{s\to1}(1-s)P_s(\Pi,\Omega)
\leq\Gamma_n^*P(\Pi,\Omega)+2\Gamma_n^*\lim_{\rho\to0^+}P(\Pi,N_\rho(\partial\Omega)),
\end{equation*}
which is $(\ref{asymptotics_nonlocal_estimate})$, once we sum the local part of the perimeter.

Their proof
relies on the fact that $\Pi$ is a polyhedron to obtain the convergence of the local part of the perimeter,
which is then used also in the estimate of the nonlocal part.  Moreover they need an approximation result to prove that
the constant is $\Gamma_n^*=\omega_{n-1}$.

They also prove, in particular
\begin{equation*}
\Gamma-\liminf_{s\to1}(1-s)P_s^L(E,\Omega)\geq\omega_{n-1}P(E,\Omega),
\end{equation*}
which is a stronger result than the first part of Theorem $\ref{Davila_conv_local}$.

\end{subsection}

\begin{subsection}{Notation and assumptions}

\begin{itemize}

\item All sets and functions considered are assumed to be Lebesgue measurable.

\item We write $A\subset\subset B$ to mean that the closure of $A$ is compact and $\overline{A}\subset B$.

\item In $\R^n$ we will usually write
$|E|=\mathcal{L}^n(E)$ for the $n$-dimensional Lebesgue measure of a set $E\subset\R^n$.

%\item By $A_h\xrightarrow{loc}A$ we mean that $\chi_{A_h}\longrightarrow\chi_A$ in $L^1_{loc}(\R^n)$,
%i.e. for every bounded open set $\Omega\subset\R^n$ we have $|(A_h\Delta A)\cap\Omega|\longrightarrow0$.

\item We write $\Ha^d$ for the $d$-dimensional Hausdorff measure, for any $d\geq0$.

\item We define the dimensional constants
\begin{equation*}
\omega_d:=\frac{\pi^\frac{d}{2}}{\Gamma\big(\frac{d}{2}+1\big)},\qquad d\geq0.
\end{equation*}
In particular, we remark that $\omega_k=\mathcal{L}^k(B_1)$ is the volume of the $k$-dimensional unit ball $B_1\subset\R^k$
and $k\,\omega_k=\Ha^{k-1}(\mathbb{S}^{k-1})$ is the surface area of the $(k-1)$-dimensional sphere
\begin{equation*}
\mathbb{S}^{k-1}=\partial B_1=\{x\in\R^k\,|\,|x|=1\}.
\end{equation*}

\item Since
\begin{equation*}
|E\Delta F|=0\quad\Longrightarrow\quad P(E,\Omega)=P(F,\Omega)\quad\textrm{and}\quad P_s(E,\Omega)=P_s(F,\Omega),
\end{equation*}
in Section 2 we implicitly identify sets up to sets of negligible Lebesgue measure.\\
%In particular, equality and inclusions of sets will usually be considered in the measure sense, e.g. $E=F$ will usually mean
%$|E\Delta F|=0$.\\
Moreover, whenever needed we can choose a particular representative for the class of $\chi_E$ in $L^1_{loc}(\R^n)$, as in the Remark below.\\
We will not make this assumption in Section 3, since the Minkowski content can be affected even by changes
in sets of measure zero, that is, in general
\begin{equation*}
|\Gamma_1\Delta\Gamma_2|=0\quad\not\Rightarrow\quad
\overline{\mathcal{M}}^r(\Gamma_1,\Omega)=\overline{\mathcal{M}}^r(\Gamma_2,\Omega)
\end{equation*}
(see Section 3 for a more detailed discussion).

\item We consider the open tubular $\rho$-neighborhood of $\partial\Omega$,
\begin{equation*}
N_\rho(\partial\Omega):=\{x\in\R^n\,|\,d(x,\partial\Omega)<\rho\}=\{|\bar{d}_\Omega|<\rho\}=\Omega_\rho\setminus\overline{\Omega_{-\rho}}
\end{equation*}
(see Appendix B).

\end{itemize}

\begin{rmk}\label{gmt_assumption}
Let $E\subset\R^n$. Up to modifying $E$ on a set of measure zero, we can assume (see Appendix C) that
\begin{equation}\label{gmt_assumption_eq}
\begin{split}
&E_1\subset E,\qquad E\cap E_0=\emptyset\\
\textrm{and}\quad\partial E=\partial^-E&=\{x\in\R^n\,|\,0<|E\cap B_r(x)|<\omega_nr^n,\,\forall\,r>0\}.
\end{split}
\end{equation}
%(inclusion and equalities holding in the usual sense, not in measure).
%This amounts to choosing a particular representative for the class of $\chi_E$ in $L^1_{loc}(\R^n)$.
%We will tacitly make this assumption in the following sections, except in Section 3, where we will explicitly talk about.
\end{rmk}

\end{subsection}

\end{section}

\begin{section}{Asymptotics as $s\to1^-$}

We say that an open set $\Omega\subset\R^n$ is an extension domain if $\exists C=C(n,s,\Omega)>0$ s.t.
for every $u\in W^{s,1}(\Omega)$ there exists $\tilde{u}\in W^{s,1}(\R^n)$ with $\tilde{u}_{|\Omega}=u$
and
\begin{equation*}
\|\tilde{u}\|_{W^{s,1}(\R^n)}\leq C\|u\|_{W^{s,1}(\Omega)}.
\end{equation*}
Every open set with bounded Lipschitz boundary is an extension domain (see \cite{HitGuide} for a proof).
For simplicity we consider $\R^n$ itself as an extension domain.

We begin with the following embedding.

\begin{prop}\label{embedding_prop}
Let $\Omega\subset\R^n$ be an extension domain. Then
$\exists C(n,s,\Omega)\geq 1$ s.t.
for every $u:\Omega\longrightarrow\R$
\begin{equation}\label{embedding_ineq}
\|u\|_{W^{s,1}(\Omega)}\leq C\|u\|_{BV(\Omega)}.
\end{equation}
In particular we have the continuous embedding
\begin{equation}
BV(\Omega)\hookrightarrow W^{s,1}(\Omega).
\end{equation}

\begin{proof}
The claim is trivially satisfied if the right hand side of $(\ref{embedding_ineq})$ is infinite, so let $u\in BV(\Omega)$.
Let $\{u_k\}\subset C^\infty(\Omega)\cap BV(\Omega)$ be an approximating sequence as in 
Theorem 1.17 of \cite{Giusti}, that is
\begin{equation*}
\|u-u_k\|_{L^1(\Omega)}\longrightarrow0\qquad\textrm{and}\qquad\lim_{k\to\infty}\int_\Omega|\nabla u_k|\,dx=|Du|(\Omega).
\end{equation*}
We only need to check that the $W^{s,1}$-seminorm of $u$ is bounded by its $BV$-norm.\\
Since $\Omega$ is an extension domain, we know (see Proposition 2.2 of \cite{HitGuide}) that
$\exists C(n,s)\geq1$ s.t.
\begin{equation*}
\|v\|_{W^{s,1}(\Omega)}\leq C\|v\|_{W^{1,1}(\Omega)}.
\end{equation*}
Then
\begin{equation*}
[u_k]_{W^{s,1}(\Omega)}\leq\|u_k\|_{W^{s,1}(\Omega)}\leq C\|u_k\|_{W^{1,1}(\Omega)}
=C\|u_k\|_{BV(\Omega)},
\end{equation*}
and hence, using Fatou's Lemma,
\begin{equation*}\begin{split}
[u]_{W^{s,1}(\Omega)}&\leq\liminf_{k\to\infty}[u_k]_{W^{s,1}(\Omega)}
\leq C\liminf_{k\to\infty}\|u_k\|_{BV(\Omega)}=C\lim_{k\to\infty}\|u_k\|_{BV(\Omega)}\\
&
=C\|u\|_{BV(\Omega)},
\end{split}\end{equation*}
proving $(\ref{embedding_ineq})$.

\end{proof}
\end{prop}

\begin{coroll}\label{embedding_fin_per_coroll}
$(i)\quad$ If $E\subset\R^n$ has finite perimeter, i.e. $\chi_E\in BV(\R^n)$, then $E$ has also finite $s$-perimeter for every
$s\in(0,1)$.\\
$(ii)\quad$ Let $\Omega\subset\R^n$ be a bounded open set with Lipschitz boundary. Then there exists $r_0>0$
s.t.
\begin{equation}\label{unif_bound_lip_frac_per}
\sup_{|r|<r_0}P_s(\Omega_r)<\infty.
\end{equation}
$(iii)\quad$ If $\Omega\subset\R^n$ is a bounded open set with Lipschitz boundary, then
\begin{equation}
P_s^{NL}(E,\Omega)\leq 2P_s(\Omega)<\infty
\end{equation}
for every $E\subset\R^n$.\\
$(iv)\quad$ Let $\Omega\subset\R^n$ be a bounded open set with Lipschitz boundary. Then
%If $E$ has finite perimeter in $\Omega$, then it has also finite $s$-perimeter in $\Omega$ for every $s\in(0,1)$.
\begin{equation}
P(E,\Omega)<\infty\qquad\Longrightarrow\qquad P_s(E,\Omega)<\infty\quad\textrm{for every }s\in(0,1).
\end{equation}

\begin{proof}
$(i)$ 
follows from
\begin{equation*}
P_s(E)=\frac{1}{2}[\chi_E]_{W^{s,1}(\R^n)}
\end{equation*}
and
previous Proposition with $\Omega=\R^n$.

$(ii)$ Let $r_0$ be as in
Proposition $\ref{bound_perimeter_unif}$
and notice that
\begin{equation*}
P(\Omega_r)=\Ha^{n-1}\big(\{\bar{d}_\Omega=r\}\big),
\end{equation*}
so that
\begin{equation*}
\|\chi_{\Omega_r}\|_{BV(\R^n)}=|\Omega_r|+\Ha^{n-1}\big(\{\bar{d}_\Omega=r\}\big).
\end{equation*}
Thus
\begin{equation*}
\sup_{|r|<r_0}P_s(\Omega_r)\leq C\Big(|\Omega_{r_0}|+\sup_{|r|<r_0}\Ha^{n-1}\big(\{\bar{d}_\Omega=r\}\big)\Big)<\infty.
\end{equation*}

$(iii)$ Notice that
\begin{equation*}\begin{split}
&\Ll_s(E\cap\Omega,\Co E\setminus\Omega)\leq \Ll_s(\Omega,\Co\Omega)=P_s(\Omega),\\
&
\Ll_s(\Co E\cap\Omega,E\setminus\Omega)\leq \Ll_s(\Omega,\Co\Omega)=P_s(\Omega),
\end{split}
\end{equation*}
and use $(\ref{unif_bound_lip_frac_per})$ (just with $\Omega_0=\Omega$).

$(iv)$ The nonlocal part of the $s$-perimeter is finite thanks to $(iii)$. As for the local part, remind that
\begin{equation*}
P(E,\Omega)=|D\chi_E|(\Omega)\qquad\textrm{and}\qquad P_s^L(E,\Omega)=\frac{1}{2}[\chi_E]_{W^{s,1}(\Omega)},
\end{equation*}
then use previous Proposition.

\end{proof}
\end{coroll}

%In this case, if we know that, we get the convergence of the (rescaled) $s$-perimeter to the classical perimeter.

\begin{subsection}{Theorem $\ref{Davila_conv_local}$, asymptotics of the local part of the $s$-perimeter}

\begin{teo}[Theorem 3' of \cite{BBM}]\label{bb}
Let $\Omega\subset\R^n$ be a smooth bounded domain. Let $u\in L^1(\Omega)$. Then
$u\in BV(\Omega)$ if and only if
\begin{equation*}
\liminf_{n\to\infty}\int_\Omega\int_\Omega\frac{|u(x)-u(y)|}{|x-y|}\rho_n(x-y)\,dxdy<\infty,
\end{equation*}
and then
\begin{equation}\label{rough}
\begin{split}
C_1|Du|(\Omega)&\leq\liminf_{n\to\infty}\int_\Omega\int_\Omega\frac{|u(x)-u(y)|}{|x-y|}\rho_n(x-y)\,dxdy\\
&
\leq\limsup_{n\to\infty}\int_\Omega\int_\Omega\frac{|u(x)-u(y)|}{|x-y|}\rho_n(x-y)\,dxdy\leq C_2|Du|(\Omega),
\end{split}
\end{equation}
for some constants $C_1$, $C_2$ depending only on $\Omega$.
\end{teo}

This result was refined by Davila

\begin{teo}[Theorem 1 of \cite{Davila} ]
Let $\Omega\subset\R^n$ be a bounded open set with Lipschitz boundary. Let $u\in BV(\Omega)$. Then
\begin{equation}\label{correct}
\lim_{k\to\infty}\int_\Omega\int_\Omega\frac{|u(x)-u(y)|}{|x-y|}\rho_k(x-y)\,dxdy=K_{1,n}|Du|(\Omega),
\end{equation}
where
\begin{equation*}
K_{1,n}=\frac{1}{n\omega_n}\int_{\mathbb{S}^{n-1}}|v\cdot e|\,d\sigma(v),
\end{equation*}
with $e\in\R^n$ any unit vector.
\end{teo}

In the above Theorems $\rho_k$ is any sequence of radial mollifiers i.e. of functions satisfying
\begin{equation}\label{rule1}
\rho_k(x)\geq0,\quad\rho_k(x)=\rho_k(|x|),\quad\int_{\R^n}\rho_k(x)\,dx=1
\end{equation}
and
\begin{equation}\label{rule2}
\lim_{k\to\infty}\int_\delta^\infty\rho_k(r)r^{n-1}dr=0\quad\textrm{for all }\delta>0.
\end{equation}

In particular, for $R$ big enough, $R>$ diam$(\Omega)$, we can consider
\begin{equation*}
\rho(x):=\chi_{[0,R]}(|x|)\frac{1}{|x|^{n-1}}
\end{equation*}
and define for any sequence $\{s_k\}\subset(0,1),\,s_k\nearrow1$,
\begin{equation*}
\rho_k(x):=(1-s_k)\rho(x)c_{s_k}\frac{1}{|x|^{s_k}},
\end{equation*}
where the $c_{s_k}$ are normalizing constants. Then
\begin{equation*}\begin{split}
\int_{\R^n}\rho_k(x)\,dx&=(1-s_k)c_{s_k}n\omega_n\int_0^R\frac{1}{r^{n-1+s_k}}r^{n-1}\,dr\\
&
=(1-s_k)c_{s_k}n\omega_n\int_0^R\frac{1}{r^{s_k}}\,dr=c_{s_k}n\omega_nR^{1-s_k},
\end{split}
\end{equation*}
and hence taking $c_{s_k}:=\frac{1}{n\omega_n}R^{s_k-1}$ gives $(\ref{rule1})$; notice that
$c_{s_k}\to\frac{1}{n\omega_n}$.\\
Also
\begin{equation*}\begin{split}
\lim_{k\to\infty}\int_\delta^\infty\rho_k(r)r^{n-1}\,dr&=
\lim_{k\to\infty}(1-s_k)c_{s_k}\int_\delta^R\frac{1}{r^{s_k}}\,dr\\
&
=\lim_{k\to\infty}c_{s_k}(R^{1-s_k}-\delta^{1-s_k})=0,
\end{split}
\end{equation*}
giving $(\ref{rule2})$.\\
With this choice we get
\begin{equation*}
\int_\Omega\int_\Omega\frac{|u(x)-u(y)|}{|x-y|}\rho_k(x-y)\,dxdy=c_{s_k}(1-s_k)[u]_{W^{s_k,1}(\Omega)}.
\end{equation*}
Then, if $u\in BV(\Omega)$, Davila's Theorem gives
\begin{equation}\label{limitperimeter}\begin{split}
\lim_{s\to1}(1-s)[u]_{W^{s,1}(\Omega)}&=\lim_{s\to1}\frac{1}{c_s}(c_s(1-s)[u]_{W^{s,1}(\Omega)})\\
&
=n\omega_nK_{1,n}|Du|(\Omega).
\end{split}
\end{equation}

\end{subsection}

\begin{subsection}{Proof of Theorem $\ref{asymptotics_teo}$}

\begin{subsubsection}{The constant $\omega_{n-1}$}

We need to compute the constant $K_{1,n}$.
Notice that we can choose $e$ in such a way that $v\cdot e=v_n$.\\
Then using spheric coordinates for $\s^{n-1}$ we obtain $|v\cdot e|=|\cos\theta_{n-1}|$
and
\begin{equation*}
d\sigma=\sin\theta_2(\sin\theta_3)^2\ldots(\sin\theta_{n-1})^{n-2}d\theta_1\ldots d\theta_{n-1},
\end{equation*}
with $\theta_1\in[0,2\pi)$ and $\theta_j\in[0,\pi)$ for $j=2,\ldots,n-1$.
Notice that
\begin{equation*}\begin{split}
\Ha^k(\s^k)&=\int_0^{2\pi}\,d\theta_1\int_0^\pi\sin\theta_2\,d\theta_2\ldots
%\int_0^\pi(\sin\theta_{k-2})^{k-3}\,d\theta_{k-2}
\int_0^\pi(\sin\theta_{k-1})^{k-2}\,d\theta_{k-1}\\
&
=\Ha^{k-1}(\s^{k-1})\int_0^\pi(\sin t)^{k-2}\,dt.
\end{split}
\end{equation*}
Then we get
\begin{equation*}
\begin{split}
\int_{\s^{n-1}}|v\cdot e|&\,d\sigma(v)=\Ha^{n-2}(\s^{n-2})\int_0^\pi(\sin t)^{n-2}|\cos t|\,dt\\
&
=\Ha^{n-2}(\s^{n-2})\Big(\int_0^\frac{\pi}{2}(\sin t)^{n-2}\cos t\,dt-\int_\frac{\pi}{2}^\pi(\sin t)^{n-2}\cos t\,dt\Big)\\
&
=\frac{\Ha^{n-2}(\s^{n-2})}{n-1}\Big(\int_0^\frac{\pi}{2}\frac{d}{dt}(\sin t)^{n-1}\,dt-\int_\frac{\pi}{2}^\pi\frac{d}{dt}(\sin t)^{n-1}\,dt\Big)\\
&
=\frac{2\Ha^{n-2}(\s^{n-2})}{n-1}.
\end{split}
\end{equation*}
Therefore
\begin{equation}
n\omega_nK_{1,n}=2\frac{\Ha^{n-2}(\s^{n-2})}{n-1}=2\Ll^{n-1}(B_1(0))=2\omega_{n-1},
\end{equation}
and hence $(\ref{limitperimeter})$ becomes
\begin{equation*}
\lim_{s\to1}(1-s)[u]_{W^{s,1}(\Omega)}=2\omega_{n-1}|Du|(\Omega),
\end{equation*}
for any $u\in BV(\Omega)$.

%Putting everything together and considering just sets, i.e. $u=\chi_E$, gives $(i)$.

\end{subsubsection}

\begin{subsubsection}{Estimating the nonlocal part of the $s$-perimeter}

We prove something slightly more general than $(\ref{asymptotics_nonlocal_estimate})$. Namely, that to estimate the nonlocal part of the $s$-perimeter
we do not necessarily need to use the sets $\Omega_\rho$: any ``regular'' approximation of $\Omega$ would do.

Let
%$\Omega\subset\R$ be a bounded open set with Lipschitz boundary. Then we can find
$A_k,\, D_k\subset\R^n$ be two sequences of
bounded open sets
with Lipschitz boundary strictly approximating $\Omega$ respectively from the inside and from the outside, that is

$(i)\quad A_k\subset A_{k+1}\subset\subset\Omega$ and $A_k\nearrow\Omega$, i.e. $\bigcup_k A_k=\Omega$,

$(ii)\quad \Omega\subset\subset D_{k+1}\subset D_k$ and $D_k\searrow\overline{\Omega}$, i.e. $\bigcap_k D_k=\overline{\Omega}$.\\
We define for every $k$
\begin{equation*}\begin{split}
&\Omega_k^+:=D_k\setminus\overline{\Omega},\qquad\Omega_k^-:=\Omega\setminus\overline{A_k}
\qquad T_k:=\Omega_k^+\cup\partial\Omega\cup\Omega_k^-,\\
&\qquad\qquad d_k:=\min\{d(A_k,\partial\Omega),\,d(D_k,\partial\Omega)\}>0.
\end{split}
\end{equation*}
In particular we can consider $\Omega_\rho$ with $\rho<0$ in place of $A_k$ and with $\rho>0$ in place of $D_k$.
Then $T_k$ would be $N_\rho(\partial\Omega)$ and $d_k=\rho$.

\begin{prop}
Let $\Omega\subset\R^n$ be a bounded open set with Lipschitz boundary and let $E\subset\R^n$
be a set having finite perimeter in $D_1$.
Then
\begin{equation}
\limsup_{s\to1}(1-s)P_s^{NL}(E,\Omega)\leq
2\omega_{n-1}\lim_{k\to\infty}P(E,T_k).
\end{equation}
In particular, if $P(E,\partial\Omega)=0$, then
\begin{equation}
\lim_{s\to1}(1-s)P_s(E,\Omega)=\omega_{n-1}P(E,\Omega).
\end{equation}

\begin{proof}
Since $\Omega$ is regular and $P(E,\Omega)<\infty$, we already know that
\begin{equation*}
\lim_{s\to1}(1-s)P_s^L(E,\Omega)=\omega_{n-1}P(E,\Omega).
\end{equation*}
Notice that, since $|D\chi_E|$ is a finite Radon measure on $D_1$ and
$T_k\searrow\partial\Omega$ as $k\nearrow\infty$, we have
\begin{equation*}
\exists\lim_{k\to\infty}P(E,T_k)=P(E,\partial\Omega).
\end{equation*}
Consider the nonlocal part of the fractional perimeter,
\begin{equation*}
P_s^{NL}(E,\Omega)=\Ll_s(E\cap\Omega,\Co E\setminus\Omega)+\Ll_s(\Co E\cap\Omega,E\setminus\Omega),
\end{equation*}
and take any $k$. Then
\begin{equation*}\begin{split}
\Ll_s(E\cap\Omega,\Co E\setminus\Omega)&=\Ll_s(E\cap\Omega,\Co E\cap\Omega_k^+)+\Ll_s(E\cap\Omega,\Co E\cap(\Co\Omega\setminus D_k))\\
&
\leq\Ll_s(E\cap\Omega,\Co E\cap\Omega_k^+)+\frac{n\omega_n}{s}|\Omega|\frac{1}{d_k^s}\\
&
\leq\Ll_s(E\cap\Omega_k^-,\Co E\cap\Omega_k^+)+2\frac{n\omega_n}{s}|\Omega|\frac{1}{d_k^s}\\
&
\leq\Ll_s(E\cap(\Omega_k^-\cup\Omega_k^+),\Co E\cap(\Omega_k^-\cup\Omega_k^+))+2\frac{n\omega_n}{s}|\Omega|\frac{1}{d_k^s}\\
&
=P^L_s(E,T_k)+2\frac{n\omega_n}{s}|\Omega|\frac{1}{d_k^s}.
\end{split}\end{equation*}
Since we can bound the other term in the same way, we get
\begin{equation}
P^{NL}_s(E,\Omega)\leq2P^L_s(E,T_k)+4\frac{n\omega_n}{s}|\Omega|\frac{1}{d_k^s}.
\end{equation}
By hypothesis we know that $T_k$ is a bounded open set with Lipschitz boundary
\begin{equation*}
\partial T_k=\partial A_k\cup\partial D_k.
\end{equation*}
Therefore using $(\ref{asymptotics_local_part})$ we have
\begin{equation*}
\lim_{s\to1}(1-s)P^L_s(E,T_k)=\omega_{n-1}P(E,T_k),
\end{equation*}
and hence
\begin{equation*}
\limsup_{s\to1}(1-s)P_s^{NL}(E,\Omega)
\leq 2\omega_{n-1}P(E,T_k).
\end{equation*}
Since this holds true for any $k$, we get the claim.

\end{proof}
\end{prop}

\end{subsubsection}

\begin{subsubsection}{Convergence in almost every $\Omega_\rho$}

Having a ``continuous'' approximating sequence (the $\Omega_\rho$) rather than numerable ones allows us to improve the previous result and obtain the second part of Theorem $\ref{asymptotics_teo}$.

We recall that De Giorgi's structure Theorem for sets of finite perimeter (see e.g. Theorem 15.9 of \cite{Maggi}) guarantees in particular that
\begin{equation*}
|D\chi_E|=\Ha^{n-1}\llcorner\partial^*E
\end{equation*}
and hence
\begin{equation*}
P(E,B)=\Ha^{n-1}(\partial^*E\cap B)\qquad\textrm{for every Borel set }B\subset\R^n,
\end{equation*}
where $\partial^*E$ is the reduced boundary of $E$.

Now suppose that $E$ has finite perimeter in $\Omega_\beta$. Then
\begin{equation*}
P(E,\partial\Omega_\delta)=\Ha^{n-1}(\partial^*E\cap\{\bar{d}_\Omega=\delta\}),
\end{equation*}
for every $\delta\in(-r_0,\beta)$.
Therefore, since
\begin{equation*}
M:=\Ha^{n-1}(\partial^*E\cap(\Omega_\beta\setminus\overline{\Omega_{-r_0}}))\leq P(E,\Omega_\beta)<\infty,
\end{equation*}
the set
\begin{equation*}
S:=\left\{\delta\in(-r_0,\beta)\,|\,P(E,\partial\Omega_\delta)>0\right\}
\end{equation*}
is at most countable.

Indeed, define
\begin{equation*}
S_k:=\Big\{\delta\in(-r_0,\beta)\,|\,\Ha^{n-1}(\partial^*E\cap\{\bar{d}_\Omega=\delta\})>\frac{1}{k}\Big\}.
\end{equation*}
Since
\begin{equation*}
\Ha^{n-1}\Big(\bigcup_{-r_0<\delta<\beta}(\partial^*E\cap\{\bar{d}_\Omega=\delta\})\Big)=M,
\end{equation*}
the number of elements in each $S_k$ is at most
\begin{equation*}
\sharp S_k\leq M\,k.
\end{equation*}
As a consequence, $S=\bigcup_k S_k$ is at most countable.\\

This concludes the proof of Theorem $\ref{asymptotics_teo}$.

\end{subsubsection}

\end{subsection}

\end{section}

\begin{section}{Irregularity of the boundary}
\begin{subsection}{The measure theoretic boundary as ``support'' of the local part of the $s$-perimeter}

First of all we show that the (local part of the) $s$-perimeter does indeed measure
a quantity related to the measure theoretic boundary.
\begin{lem}
Let $E\subset\R^n$ be a set of locally finite $s$-perimeter. Then
\begin{equation}
\partial^-E=\{x\in\R^n\,|\,P_s^L(E,B_r(x))>0\textrm{ for every }r>0\}.
\end{equation}
\begin{proof}
The claim follows from the following observation. Let $A,\,B\subset\R^n$ s.t. $A\cap B=\emptyset$; then
\begin{equation*}
\Ll_s(A,B)=0\quad\Longleftrightarrow\quad|A|=0\quad\textrm{or}\quad|B|=0.
\end{equation*}
Therefore
\begin{equation*}\begin{split}
x\in\partial^-E&\quad\Longleftrightarrow\quad
|E\cap B_r(x)|>0\textrm{  and }|\Co E\cap B_r(x)|>0\quad\forall\,r>0\\
&
\quad\Longleftrightarrow\quad
\Ll_s(E\cap B_r(x),\Co E\cap B_r(x))>0\quad\forall\,r>0.
\end{split}\end{equation*}

\end{proof}
\end{lem}

This characterization of $\partial^-E$ can be thought of as a fractional analogue
of $(\ref{support_perimeter})$. However we can not really think of $\partial^-E$ as the support of
\begin{equation*}
P_s^L(E,-):\Omega\longmapsto P_s^L(E,\Omega),
\end{equation*}
in the sense that, in general
\begin{equation*}
\partial^-E\cap\Omega=\emptyset\quad\not\Rightarrow\quad P_s^L(E,\Omega)=0.
\end{equation*}
For example, consider $E:=\{x_n\leq0\}\subset\R^n$ and notice that $\partial^-E=\{x_n=0\}$.
Let $\Omega:=B_1(2e_n)\cup B_1(-2e_n)$. Then $\partial^-E\cap\Omega=\emptyset$, but
\begin{equation*}
P_s^L(E,\Omega)=\Ll_s(B_1(2e_n),B_1(-2e_n))>0.
\end{equation*}

On the other hand, the only obstacle is the non connectedness of the set $\Omega$ and
indeed we obtain the following
\begin{prop}
Let $E\subset\R^n$ be a set of locally finite $s$-perimeter and let $\Omega\subset\R^n$ be an open set.
Then
\begin{equation*}
\partial^-E\cap\Omega\not=\emptyset\quad\Longrightarrow\quad P_s^L(E,\Omega)>0.
\end{equation*}
Moreover, if $\Omega$ is connected
\begin{equation*}
\partial^-E\cap\Omega=\emptyset\quad\Longrightarrow\quad P_s^L(E,\Omega)=0.
\end{equation*}
Therefore, if $\widehat{\mathcal O}(\R^n)$ denotes the family of bounded and connected open sets, then
$\partial^-E$ is the ``support'' of
\begin{equation*}\begin{split}
P_s^L(E,-):\,&\widehat{\mathcal O}(\R^n)\longrightarrow [0,\infty)\\
&
\Omega\longmapsto P_s^L(E,\Omega),
\end{split}\end{equation*}
in the sense that, if $\Omega\in\widehat{\mathcal O}(\R^n)$, then
\begin{equation*}
P_s^L(E,\Omega)>0\quad\Longleftrightarrow\quad\partial^-E\cap\Omega\not=\emptyset.
\end{equation*}
\begin{proof}
Let $x\in\partial^-E\cap\Omega$. Since $\Omega$ is open, we have $B_r(x)\subset\Omega$ for some $r>0$
and hence
\begin{equation*}
P_s^L(E,\Omega)\geq P_s^L(E,B_r(x))>0.
\end{equation*}
Let $\Omega$ be connected and suppose $\partial^-E\cap\Omega=\emptyset$.
We have the partition of $\R^n$ as $\R^n=E_0\cup\partial^-E\cup E_1$ (see Appendix C). Thus we
can write $\Omega$ as the disjoint union
\begin{equation*}
\Omega=(E_0\cap\Omega)\cup(E_1\cap\Omega).
\end{equation*}
However, since $\Omega$ is connected and both $E_0$ and $E_1$ are open,
we must have $E_0\cap\Omega=\emptyset$ or $E_1\cap\Omega=\emptyset$.
Now, if $E_0\cap\Omega=\emptyset$ (the other case is analogous),
then $\Omega\subset E_1$ and hence $|\Co E\cap\Omega|=0$. Thus
\begin{equation*}
P_s^L(E,\Omega)=\Ll_s(E\cap\Omega,\Co E\cap\Omega)=0.
\end{equation*}

\end{proof}
\end{prop}

\end{subsection}

\begin{subsection}{A notion of fractal dimension}

Let $\Omega\subset\R^n$ be an open set. Then
\begin{equation*}
t>s\qquad\Longrightarrow\qquad W^{t,1}(\Omega)\hookrightarrow W^{s,1}(\Omega),
\end{equation*}
(see e.g. Proposition 2.1 of \cite{HitGuide}). As a consequence, for every $u:\Omega\longrightarrow\R$ there exists a unique
$R(u)\in[0,1]$ s.t.
\begin{equation*}
[u]_{W^{s,1}(\Omega)}\quad\left\{\begin{array}{cc}
<\infty,& \forall\,s\in(0,R(u))\\
=\infty, &\forall\,s\in(R(u),1)
\end{array}\right.
\end{equation*}
that is
\begin{equation}\begin{split}\label{frac_range}
R(u)&=\sup\left\{s\in(0,1)\,\big|\,[u]_{W^{s,1}(\Omega)}<\infty\right\}\\
&
=\inf\left\{s\in(0,1)\,\big|\,[u]_{W^{s,1}(\Omega)}=\infty\right\}.
\end{split}
\end{equation}

In particular, exploiting this result for characteristic functions, in \cite{Visintin} the author suggested the following definition of fractal dimension.
\begin{defin}
Let $\Omega\subset\R^n$ be an open set and let $E\subset\R^n$. If $\partial^- E\cap\Omega\not=\emptyset$, we define
\begin{equation}
\Dim_F(\partial^- E,\Omega):=n-R(\chi_E),
\end{equation}
the fractal dimension of $\partial^- E$ in $\Omega$, relative to the fractional perimeter.\\
If $\Omega=\R^n$, we drop it in the formulas.
\end{defin}

Notice that in the case of sets $(\ref{frac_range})$ becomes
\begin{equation}
\begin{split}\label{frac_range_sets}
R(\chi_E)&=\sup\left\{s\in(0,1)\,\big|\,P_s^L(E,\Omega)<\infty\right\}\\
&
=\inf\left\{s\in(0,1)\,\big|\,P_s^L(E,\Omega)=\infty\right\}.
\end{split}
\end{equation}

In particular we can take $\Omega$ to be the whole of $\R^n$, or a bounded open set with Lipschitz boundary.\\
In the first case the local part of the fractional perimeter coincides with the whole fractional perimeter, while in the second case we know that we can bound the nonlocal part with $2P_s(\Omega)<\infty$ for every $s\in(0,1)$. Therefore in both cases in $(\ref{frac_range_sets})$
we can as well take the whole fractional perimeter $P_s(E,\Omega)$ instead of just the local part.\\

Now we give a proof of the relation $(\ref{intro_dim_ineq})$ (obtained in \cite{Visintin}).\\
For simplicity, given $\Gamma\subset\R^n$ we set
\begin{equation}\label{neigh_mink_def}
\bar{N}_\rho^\Omega(\Gamma):=\overline{N_\rho(\Gamma)}\cap\Omega
=\{x\in\Omega\,|\,d(x,\Gamma)\leq\rho\},
\end{equation}
for any $\rho>0$.

\begin{prop}\label{vis_prop}
Let $\Omega\subset\R^n$ be a bounded open set. Then for every $E\subset\R^n$ s.t. $\partial^- E\cap\Omega\not=\emptyset$ and $\overline{\Dim}_\mathcal{M}(\partial^-E,\Omega)\geq n-1$ we have
\begin{equation}
\Dim_F(\partial^-E,\Omega)\leq\overline{\Dim}_\mathcal{M}(\partial^-E,\Omega).
\end{equation}
\begin{proof}
By hypothesis we have
\begin{equation*}
\overline{\Dim}_\mathcal{M}(\partial^-E,\Omega)=n-\inf\big\{r\in(0,1)\,|\,\overline{\mathcal{M}}^{n-r}(\partial^-E,\Omega)=\infty\big\},
\end{equation*}
and we need to show that
\begin{equation*}
\inf\big\{r\in(0,1)\,|\,\overline{\mathcal{M}}^{n-r}(\partial^-E,\Omega)=\infty\big\}
\leq
\sup\{s\in(0,1)\,|\,P_s^L(E,\Omega)<\infty\}.
\end{equation*}
Up to modifying $E$ on a set of Lebesgue measure zero
we can suppose that $\partial E=\partial^-E$, as in Remark $\ref{gmt_assumption}$. Notice that this does not affect the
$s$-perimeter.

Now for any $s\in(0,1)$
\begin{equation*}\begin{split}
2P_s^L(E,\Omega)&=\int_\Omega\,dx\int_\Omega\frac{|\chi_E(x)-\chi_E(y)|}{|x-y|^{n+s}}\,dy\\
&
=\int_\Omega dx\int_0^\infty d\rho\int_{\partial B_\rho(x)\cap\Omega}\frac{|\chi_E(x)-\chi_E(y)|}{|x-y|^{n+s}}\,d\Ha^{n-1}(y)\\
&
=\int_\Omega dx\int_0^\infty\frac{d\rho}{\rho^{n+s}}\int_{\partial B_\rho(x)\cap\Omega}|\chi_E(x)-\chi_E(y)|\,d\Ha^{n-1}(y).
\end{split}
\end{equation*}
Notice that
\begin{equation*}
d(x,\partial E)>\rho\quad\Longrightarrow\quad\chi_E(y)=\chi_E(x),\quad\forall\,y\in\overline{B_\rho(x)},
\end{equation*}
and hence
\begin{equation*}\begin{split}
\int_{\partial B_\rho(x)\cap\Omega}|\chi_E(x)-\chi_E(y)|\,d\Ha^{n-1}(y)&
\leq\int_{\partial B_\rho(x)\cap\Omega}\chi_{\bar{N}_\rho(\partial E)}(x)\,d\Ha^{n-1}(y)\\
&
\leq n\omega_n\rho^{n-1}\chi_{\bar{N}_\rho(\partial E)}(x).
\end{split}
\end{equation*}
Therefore
\begin{equation}\label{visintin_pf}
2P_s^L(E,\Omega)\leq n\omega_n\int_0^\infty\frac{d\rho}{\rho^{1+s}}\int_\Omega
\chi_{\bar{N}_\rho(\partial E)}(x)
=n\omega_n\int_0^\infty\frac{|\bar{N}^\Omega_\rho(\partial E)|}{\rho^{1+s}}\,d\rho.
\end{equation}
We prove the following\\
CLAIM
\begin{equation}\label{visintin_proof}
\overline{\mathcal{M}}^{n-r}(\partial E,\Omega)<\infty\quad\Longrightarrow\quad P_s^L(E,\Omega)<\infty,\quad\forall\,s\in(0,r).
\end{equation}
Indeed
\begin{equation*}
\limsup_{\rho\to0}\frac{|\bar{N}^\Omega_\rho(\partial E)|}{\rho^r}<\infty\quad\Longrightarrow\quad\exists\,C>0\textrm{ s.t. }
\sup_{\rho\in(0,C]}\frac{|\bar{N}^\Omega_\rho(\partial E)|}{\rho^r}\leq M<\infty.
\end{equation*}
Then
\begin{equation*}\begin{split}
2P_s^L(E,\Omega)&\leq n\omega_n\Big\{\int_0^C\frac{|\bar{N}^\Omega_\rho(\partial E)|}{\rho^{1-(r-s)+r}}\,d\rho
+\int_C^\infty\frac{|\bar{N}^\Omega_\rho(\partial E)|}{\rho^{1+s}}\,d\rho\Big\}\\
&
\leq n\omega_n\Big\{
M\int_0^C\frac{1}{\rho^{1-(r-s)}}\,d\rho+|\Omega|\int_C^\infty\frac{1}{\rho^{1+s}}\,d\rho
\Big\}\\
&
=n\omega_n\Big\{
\frac{M}{r-s}C^{r-s}+\frac{|\Omega|}{sC^s}
\Big\}<\infty,
\end{split}\end{equation*}
proving the claim.\\
This implies
\begin{equation*}
r\leq\sup\{s\in(0,1)\,|\,P_s^L(E,\Omega)<\infty\},
\end{equation*}
for every $r\in(0,1)$ s.t. $\overline{\mathcal{M}}^{n-r}(\partial E,\Omega)<\infty$.\\
Thus for $\epsilon>0$ very small, we have
\begin{equation*}
\inf\big\{r\in(0,1)\,|\,\overline{\mathcal{M}}^{n-r}(\partial^-E,\Omega)=\infty\big\}-\epsilon
\leq\sup\{s\in(0,1)\,|\,P_s^L(E,\Omega)<\infty\}.
\end{equation*}
Letting $\epsilon$ tend to zero, we conclude the proof.

\end{proof}
\end{prop}

In particular, if $\Omega$ has Lipschitz boundary we obtain
\begin{coroll}\label{fractal_dim_coroll}
Let $\Omega\subset\R^n$ be a bounded open set with Lipschitz boundary. Let $E\subset\R^n$ s.t. $\partial^-E\cap\Omega\not=\emptyset$ and
$\overline{\Dim}_\mathcal{M}(\partial^-E,\Omega)\in[n-1,n)$. Then
\begin{equation}\label{fractal_per}
P_s(E,\Omega)<\infty\qquad\textrm{for every }s\in\left(0,n-\overline{\Dim}_\mathcal{M}(\partial^-E,\Omega)\right).
\end{equation}
\end{coroll}

\begin{rmk}\label{fractal_dim_rmk}
Actually, previous Proposition and Corollary still work when $\Omega=\R^n$, provided the set $E$ we are considering is bounded.\\
Indeed, if $E$ is bounded, we can apply previous results with $\Omega=B_R$ s.t. $E\subset\Omega$. Moreover, since $\Omega$
has a regular boundary, as remarked above we can take the whole $s$-perimeter in
$(\ref{frac_range_sets})$, instead of just the local part. But then, since $P_s(E,\Omega)=P_s(E)$, we see that
\begin{equation*}
\Dim_F(\partial^-E,\Omega)=\Dim_F(\partial^-E,\R^n).
\end{equation*}
\end{rmk}

\begin{subsubsection}{The measure theoretic boundary of a set of locally finite $s$-perimeter (in general) is not rectifiable}

These results show that a set $E$ can have finite fractional perimeter even if its boundary is really irregular, unlike what happens with a Caccioppoli set and its reduced boundary, which is locally $(n-1)$-rectifiable.\\

\noindent
Indeed, if $\Omega\subset\R^n$ is a bounded open set with Lipschitz boundary and $E\subset\R^n$ is s.t.
$\emptyset\not=\partial^-E\cap\Omega$ is not $(n-1)$-rectifiable, with
$\overline{\Dim}_\mathcal{M}(\partial^-E,\Omega)\in(n-1,n)$,
thanks to previous Corollary we have $P_s(E,\Omega)<\infty$ for every $s\in(0,\sigma)$.

We give some examples of this kind of sets in the following Sections.\\
In particular, the von Koch snowflake $S\subset\R^2$ has finite $s$-perimeter for every $s\in(0,\sigma)$, but $\partial^-S=\partial S$
is not locally $(n-1)$-rectifiable.\\
Actually, because of the self-similarity of the von Koch curve, there is no part of $\partial S$ which is rectifiable (see below).\\

On the other hand, De Giorgi's structure Theorem (see e.g. Theorem 15.9 and Corollary 16.1 of \cite{Maggi})
says that if a set $E\subset\R^n$ has locally finite perimeter, then its reduced boundary $\partial^*E$ is
locally $(n-1)$-rectifiable.\\
Moreover the reduced boundary is dense in the measure theoretic boundary, which is the support of
the Radon measure $|D\chi_E|$,
\begin{equation*}
\overline{\partial^*E}=\partial^-E=\textrm{supp }|D\chi_E|.
\end{equation*}

This underlines a deep difference between the classical perimeter and the $s$-perimeter, which can indeed be thought of as a
``fractional'' perimeter.\\
Namely, having (locally) finite classical perimeter implies the regularity of an ``important'' portion of
the (measure theoretic) boundary. On the other hand, a set can have a fractal, nowhere rectifiable boundary
and still have (locally) finite $s$-perimeter.

\end{subsubsection}

\begin{subsubsection}{Remarks about the Minkowski content of $\partial^-E$}

In the beginning of the proof of Proposition $\ref{vis_prop}$
we chose a particular representative for the class of $E$ in order to have $\partial E=\partial^-E$.
This can be done since it does not affect the $s$-perimeter
and we are already considering the Minkowski dimension of $\partial^-E$.

On the other hand, if we consider a set $F$ s.t. $|E\Delta F|=0$,
we can use the same proof to obtain the inequality
\begin{equation*}
\Dim_F(\partial^-E,\Omega)\leq\overline{\Dim}_\mathcal{M}(\partial F,\Omega).
\end{equation*}
It is then natural to ask whether we can find a ``better'' representative $F$,
whose (topological) boundary $\partial F$ has Minkowski dimension strictly smaller than that
of $\partial^-E$.

First of all, we remark that the Minkowski content can be influenced by changes in sets of measure zero.
Roughly speaking, this is because the Minkowski content is not a purely measure theoretic notion,
but rather a combination of metric and measure.

For example, let $\Gamma\subset\R^n$ and define $\Gamma':=\Gamma\cup\mathbb Q^n$.
Then $|\Gamma\Delta\Gamma'|=0$, but
$N_\delta(\Gamma')=\R^n$ for every $\delta>0$.

In particular, considering different representatives for $E$
we will get different topological boundaries and hence different Minkowski dimensions.

However, since the measure theoretic boundary minimizes the size of the topological boundary, that is
\begin{equation*}
\partial^-E=\bigcap_{|F\Delta E|=0}\partial F,
\end{equation*}
(see Appendix C), it minimizes also the Minkowski dimension.\\
Indeed, for every $F$ s.t. $|F\Delta E|=0$ we have
\begin{equation*}\begin{split}
\partial^-E\subset\partial F&\quad\Longrightarrow\quad \bar N_\rho^\Omega(\partial^-E)\subset
\bar N_\rho^\Omega(\partial F)\\
&
\quad\Longrightarrow\quad
\overline{\mathcal M}^r(\partial^-E,\Omega)\leq
\overline{\mathcal M}^r(\partial F,\Omega)\\
&
\quad\Longrightarrow\quad
\overline{\Dim}_\mathcal{M}(\partial^-E,\Omega)\leq
\overline{\Dim}_\mathcal{M}(\partial F,\Omega).
\end{split}
\end{equation*}

\end{subsubsection}

\end{subsection}

\begin{subsection}{Fractal dimension of the von Koch snowflake}

The von Koch snowflake $S\subset\R^2$ is an example of bounded open set with fractal boundary,
for which the Minkowski dimension and the fractal dimension introduced above coincide.

Moreover its boundary is ``nowhere rectifiable'', in the sense that
$\partial S\cap B_r(p)$ is not $(n-1)$-rectifiable for any $r>0$
and $p\in\partial S$.\\

First of all we construct the von Koch curve. Then the snowflake is made of three
von Koch curves.

Let $\Gamma_0$ be a line segment of unit length. 
The set $\Gamma_1$ consists of the four segments obtained by removing the middle third of $\Gamma_0$ and replacing it by the other two sides of the equilateral triangle based on the removed segment.\\
We construct $\Gamma_2$ by applying the same procedure to each of the segments in $\Gamma_1$ and so on.
Thus $\Gamma_k$ comes from replacing the middle third of each straight line segment of $\Gamma_{k-1}$ by the other two sides of an equilateral triangle.

As $k$ tends to infinity, the sequence of polygonal curves $\Gamma_k$ approaches a limiting curve $\Gamma$, called the von Koch curve.\\
If we start with an equilateral triangle with unit length side and perform the same construction on all three sides, we obtain the von Koch snowflake $\Sigma$.\\
Let $S$ be the bounded region enclosed by $\Sigma$, so that $S$ is open and $\partial S=\Sigma$. We still
call $S$ the von Koch snowflake.\\

Now we calculate the (Minkowski) dimension of $\Gamma$ using the box-counting dimensions (see Appendix D).\\
The idea is to exploit the self-similarity of $\Gamma$ and consider covers made of squares with side $\delta_k=3^{-k}$.

The key observation is that $\Gamma$ can be covered by three squares of length $1/3$ (and cannot be covered by only two),
so that $\mathcal{N}(\Gamma,1/3)=3$.\\
Then consider $\Gamma_1$. We can think of $\Gamma$ as being made of four von Koch curves starting from the set $\Gamma_1$ and with initial segments of length $1/3$ instead of 1. Therefore we can cover each of these four pieces with three squares of side $1/9$, so that $\Gamma$ can be covered with $3\cdot4$ squares of length $1/9$ (and not one less) and $\mathcal{N}(\Gamma,1/9)=4\cdot3$.

We can repeat the same argument starting from $\Gamma_2$ to get $\mathcal{N}(\Gamma,1/27)=4^2\cdot3$, and so on.
In general we obtain
\begin{equation*}
\mathcal{N}(\Gamma,3^{-k})=4^{k-1}\cdot3.
\end{equation*}
Then, taking logarithms we get
\begin{equation*}
\frac{\log\mathcal{N}(\Gamma,3^{-k})}{-\log3^{-k}}=\frac{\log3+(k-1)\log4}{k\log3}\longrightarrow\frac{\log4}{\log3},
\end{equation*}
so that $\Dim_\mathcal{M}(\Gamma)=\frac{\log4}{\log3}$.

Notice that the Minkowski dimensions of the snowflake and of the curve are the same.
Moreover it can be shown that the Hausdorff dimension of the von Koch curve is equal to its Minkowski dimension, so we
obtain
\begin{equation}\label{dime_Koch_snow}
\Dim_\Ha(\Sigma)=\Dim_\mathcal M(\Sigma)=\frac{\log4}{\log3}
\end{equation}

Now we explain how to construct $S$ in a recursive way and we prove that
\begin{equation*}
\partial^-S=\partial S=\Sigma.
\end{equation*}

As starting point for the snowflake take the
equilateral triangle $T$ of side 1, with baricenter in the origin and a vertex on the $y$-axis, $P=(0,t)$ with $t>0$.\\
Then $T_1$ is made of three triangles of side $1/3$, $T_2$ of $3\cdot4$ triangles of side $1/3^2$ and so on.\\
In general $T_k$ is made of $3\cdot4^{k-1}$ triangles of side $1/3^k$, call them $T_k^1,\ldots,T_k^{3\cdot4^{k-1}}$.
Let $x^i_k$ be the baricenter of $T_k^i$ and $P_k^i$ the vertex which does not touch $T_{k-1}$.

Then $S=T\cup\bigcup T_k$. Also notice that $T_k$ and $T_{k-1}$ touch only on a set of measure zero.

For each triangle
$T^i_k$ there exists a rotation $\mathcal{R}_k^i\in SO(n)$ s.t.
\begin{equation*}
T_k^i=F_k^i(T):=\mathcal{R}_k^i\Big(\frac{1}{3^k}T\Big)+x_k^i.
\end{equation*}
We choose the rotations so that $F_k^i(P)=P_k^i$.

Notice that for each triangle $T_k^i$ we can find a small ball which is contained in the complementary
of the snowflake,
$B_k^i\subset\Co S$,
and touches the triangle in the vertex $P_k^i$. Actually these balls can be obtained as the images
of the affine transformations $F_k^i$
of a fixed
ball $B$.

To be more precise, fix a small ball contained in the
complementary of $T$, which has the center on the $y$-axis
and touches $T$ in the vertex $P$, say $B:=B_{1/1000}(0,t+1/1000)$. Then
\begin{equation}\label{koch3}
%B_{3^{-k}}(x+x_k^i)=\mathcal{R}_k^i\Big(\frac{1}{3^k}B_1(x)\Big)+x_k^i\subset\Co S,
B_k^i:=F_k^i(B)\subset\Co S
\end{equation}
for every $i,\,k$. To see this, imagine constructing the snowflake $S$ using the same affine transformations $F_k^i$
but starting with $T\cup B$ in place of $T$.\\

We know that $\partial^-S\subset\partial S$ (see Appendix C).\\
On the other hand, let $p\in\partial S$. Then
every ball $B_\delta(p)$ contains at least a triangle $T^i_k\subset S$ and its corresponding ball
$B^i_k\subset\Co S$ (and actually infinitely many). Therefore $0<|B_\delta(p)\cap S|<\omega_n\delta^n$
for every $\delta>0$ and hence $p\in\partial^-S$.

\begin{proof}[proof of Theorem $\ref{von_koch_snow}$]
Since $S$ is bounded, its boundary is $\partial^-S=\Sigma$, and
$\Dim_\mathcal M(\Sigma)=\frac{\log4}{\log3}$,
we obtain $(\ref{koch1})$ from Corollary $\ref{fractal_dim_coroll}$
and Remark $\ref{fractal_dim_rmk}$.

Exploiting the construction of $S$ given above and $(\ref{koch3})$ we prove $(\ref{koch2})$.\\
We have
\begin{equation*}\begin{split}
P_s(S)&=\Ll_s(S,\Co S)=\Ll_s(T,\Co S)+\sum_{k=1}^\infty\Ll_s(T_k,\Co S)\\
&
=\Ll_s(T,\Co S)+\sum_{k=1}^\infty\sum_{i=1}^{3\cdot4^{k-1}}\Ll_s(T_k^i,\Co S)
\geq\sum_{k=1}^\infty\sum_{i=1}^{3\cdot4^{k-1}}\Ll_s(T_k^i,\Co S)\\
&
\geq\sum_{k=1}^\infty\sum_{i=1}^{3\cdot4^{k-1}}\Ll_s(T_k^i,B_k^i)\qquad\textrm{(by }(\ref{koch3}))\\
&
=\sum_{k=1}^\infty\sum_{i=1}^{3\cdot4^{k-1}}\Ll_s(F_k^i(T),F_k^i(B))\\
&
=\sum_{k=1}^\infty\sum_{i=1}^{3\cdot4^{k-1}}\Big(\frac{1}{3^k}\Big)^{2-s}\Ll_s(T,B)\qquad\textrm{(by Proposition }\ref{elementary_properties})\\
&
=\frac{3}{3^{2-s}}\Ll_s(T,B)\sum_{k=0}^\infty\Big(\frac{4}{3^{2-s}}\Big)^k.
\end{split}\end{equation*}
We remark that
\begin{equation*}
\Ll_s(T,B)\leq\Ll_s(T,\Co T)=P_s(T)<\infty,
\end{equation*}
for every $s\in(0,1)$.

To conclude, notice that the last series is divergent if $s\geq2-\frac{\log4}{\log3}$.

\end{proof}

Exploiting the self-similarity of the von Koch curve, we show that the fractal dimension
of $S$ is the same in every open set which contains a point of $\partial S$.

%We exploit this fact to prove (the well known property) that $\partial S$ is nowhere rectifiable.

\begin{coroll}\label{koch_coroll}
Let $S\subset\R^2$ be the von Koch snowflake. Then
\begin{equation*}
\Dim_F(\partial S,\Omega)=\frac{\log4}{\log3}
\end{equation*}
for every open set $\Omega$ s.t. $\partial S\cap\Omega\not=\emptyset$.
%Therefore, if $p\in\partial S$ and $r>0$, then $\partial S\cap B_r(p)$ is not $(n-1)$-rectifiable.

\begin{proof}
Since $P_s(S,\Omega)\leq P_s(S)$, we have
\begin{equation*}
P_s(S,\Omega)<\infty,\qquad\forall\,s\in\Big(0,2-\frac{\log4}{\log3}\Big).
\end{equation*}
On the other hand, if $p\in\partial S\cap\Omega$, then $B_r(p)\subset\Omega$
for some $r>0$.
Now notice that $B_r(p)$ contains a rescaled version of the von Koch curve, including
all the triangles $T_k^i$ which constitute it and the relative
balls $B_k^i$. We can thus repeat the argument above to obtain
\begin{equation*}
P_s(S,\Omega)\geq P_s(S,B_r(p))=\infty,\qquad\forall\,s\in\Big[2-\frac{\log4}{\log3},1\Big).
\end{equation*}

\end{proof}
\end{coroll}

\end{subsection}

\begin{subsection}{Self-similar fractal boundaries}

The von Koch curve is a well known example of a family of rather ``regular'' fractal sets,
the self-similar fractal sets (see e.g. Section 9 of \cite{Falconer} for the proper definition and the main properties).

Many examples of this kind of sets can be constucted in a recursive way similar to that of the von Koch snowflake.

To be more precise, 
we start with a bounded open set $T_0\subset\R^n$ with finite perimeter $P(T_0)<\infty$, which is, roughly speaking, our basic ``building block''.

Then we go on inductively by adding roto-translations of a scaling of the building block $T_0$,
i.e. sets of the form
\begin{equation*}
T_k^i=F_k^i(T_0):=\mathcal{R}_k^i\big(\lambda^{-k}T_0\big)+x_k^i,
\end{equation*}
where $\lambda>1$, $k\in\mathbb N$, $1\leq i\leq ab^{k-1}$, with $a,\,b\in\mathbb N$,
$\mathcal{R}_k^i\in SO(n)$ and $x_k^i\in\R^n$.
We ask that these sets do not overlap, i.e.
\begin{equation*}
|T^i_k\cap T^j_h|=0,\qquad\textrm{if }i\not=j.
\end{equation*}
Then we define
\begin{equation}\label{frac_ind_def}
T_k:=\bigcup_{i=1}^{ab^{k-1}}T_k^i\qquad\textrm{and}\qquad T:=\bigcup_{k=1}^\infty T_k.
\end{equation}
The final set $E$ is either
\begin{equation}
E:=T_0\cup\bigcup_{k\geq1}\bigcup_{i=1}^{ab^{k-1}}T^i_k,\quad\textrm{or}\quad
E:=T_0\setminus\Big(\bigcup_{k\geq1}\bigcup_{i=1}^{ab^{k-1}}T^i_k\Big).
\end{equation}
For example, the von-Koch snowflake is obtained by adding pieces.

Examples obtained by removing the $T_k^i$'s
are the middle Cantor set $E\subset\R$, the Sierpinski triangle $E\subset\R^2$
and Menger sponge $E\subset\R^3$.\\

We will consider just the set $T$ and exploit the same argument used for the von Koch snowflake
to compute the fractal dimension related to the $s$-perimeter.\\
However, the Cantor set, the Sierpinski triangle and the Menger sponge are s.t. $|E|=0$, i.e. $|T_0\Delta T|=0$.\\
Therefore both the perimeter and the $s$-perimeter do not notice the fractal nature of the (topological) boundary of $T$
and indeed, since $P(T)=P(T_0)<\infty$, we get $P_s(T)<\infty$ for every $s\in(0,1)$.
For example, in the case of the Sierpinski triangle, $T_0$ is an equilateral triangle
and $\partial^-T=\partial T_0$, even if $\partial T$ is a self-similar fractal.

Roughly speaking, the problem in these cases is that there is not room enough to find a small ball $B_k^i=F_k^i(B)\subset\Co T$
near each piece $T_k^i$.

Therefore, we will make the additional assumption that
\begin{equation}\label{add_frac_self_hp}
\exists\,S_0\subset\Co T\quad\textrm{s.t. }|S_0|>0\quad\textrm{and }S_k^i:=F_k^i(S_0)\subset\Co T\quad\forall\,k,\,i.
\end{equation}
We remark that it is not necessary to ask that these sets do not overlap.
%even if this will be the case in the following examples. %we give below.\\
%It might seem that this condition is too restrictive to be satisfied by interesting sets, since 
%Even if the examples above (except the snowflake) do not satisfy this condition,
%but this is not the case.

%It is true that many examples of self-similar fractals do not satisfy $(\ref{add_frac_self_hp})$.
%However b
Below we give some examples on how to construct %(in an easy way) 
sets which satisfy this additional hypothesis starting with sets which do not,
like the Sierpinski triangle,
without altering their ``structure''.

\begin{teo}\label{fractal_bdary_selfsim_dim}
Let $T\subset\R^n$ be a set which can be written as in $(\ref{frac_ind_def})$.
If $\frac{\log b}{\log\lambda}\in(n-1,n)$ and $(\ref{add_frac_self_hp})$ holds true, then
\begin{equation*}
P_s(T)<\infty,\qquad\forall\,s\in\Big(0,n-\frac{\log b}{\log\lambda}\Big)
\end{equation*}
and
\begin{equation*}
P_s(T)=\infty,\qquad\forall\,s\in\Big[n-\frac{\log b}{\log\lambda},1\Big).
\end{equation*}
Thus
\begin{equation}
\Dim_F(\partial^-T)=\frac{\log b}{\log\lambda}.
\end{equation}

\begin{proof}
Arguing as we did with the von Koch snowflake, we show that $P_s(T)$ is
bounded both from above and from below by the series
\begin{equation*}
\sum_{k=0}^\infty\Big(\frac{b}{\lambda^{n-s}}\Big)^k,
\end{equation*}
which converges if and only if $s<n-\frac{\log b}{\log\lambda}$.

Indeed
\begin{equation*}\begin{split}
P_s(T)&=\Ll_s(T,\Co T)=\sum_{k=1}^\infty\sum_{i=1}^{ab^{k-1}}\Ll_s(T_k^i,\Co T)\\
&
\leq
\sum_{k=1}^\infty\sum_{i=1}^{ab^{k-1}}\Ll_s(T_k^i,\Co T_k^i)
=
\sum_{k=1}^\infty\sum_{i=1}^{ab^{k-1}}\Ll_s(F_k^i(T_0),F_k^i(\Co T_0))\\
&
=\frac{a}{\lambda^{n-s}}\Ll_s(T_0,\Co T_0)\sum_{k=0}^\infty\Big(\frac{b}{\lambda^{n-s}}\Big)^k,
\end{split}\end{equation*}
and
\begin{equation*}\begin{split}
P_s(T)&=\Ll_s(T,\Co T)=\sum_{k=1}^\infty\sum_{i=1}^{ab^{k-1}}\Ll_s(T_k^i,\Co T)\\
&
\geq
\sum_{k=1}^\infty\sum_{i=1}^{ab^{k-1}}\Ll_s(T_k^i,S_k^i)
=
\sum_{k=1}^\infty\sum_{i=1}^{ab^{k-1}}\Ll_s(F_k^i(T_0),F_k^i(S_0))\\
&
=\frac{a}{\lambda^{n-s}}\Ll_s(T_0,S_0)\sum_{k=0}^\infty\Big(\frac{b}{\lambda^{n-s}}\Big)^k.
\end{split}\end{equation*}
Also notice that, since $P(T_0)<\infty$, we have
\begin{equation*}
\Ll_s(T_0,S_0)\leq\Ll_s(T_0,\Co T_0)=P_s(T_0)<\infty,
\end{equation*}
for every $s\in(0,1)$.

\end{proof}
\end{teo}

Now suppose that $T$ does not satisfy $(\ref{add_frac_self_hp})$.
Then we can obtain a set $T'$ which does, simply by removing
a part $S_0$ of the building block $T_0$.\\
To be more precise, let $S_0\subset T_0$ be s.t. $|S_0|>0$, $|T_0\setminus S_0|>0$ and $P(T_0\setminus S_0)<\infty$.
Then define a new
building block $T'_0:=T_0\setminus S_0$ and the set
\begin{equation*}
T':=\bigcup_{k=1}^\infty\bigcup_{i=1}^{ab^{k-1}}F_k^i(T'_0).
\end{equation*}
This new set has exactly the same structure of $T$, since we are using the same
collection $\{F_k^i\}$ of affine maps. 

Notice that
\begin{equation*}
S_0\subset T_0\quad\Longrightarrow\quad F_k^i(S_0)\subset F_k^i(T_0),
\end{equation*}
and
\begin{equation*}
F_k^i(T'_0)=F_k^i(T_0)\setminus F_k^i(S_0),
\end{equation*}
for every $k,\,i$.
Thus
\begin{equation*}
T'=T\setminus\Big(\bigcup_{k=1}^\infty\bigcup_{i=1}^{ab^{k-1}}F_k^i(S_0)\Big)
\end{equation*}
satisfies $(\ref{add_frac_self_hp})$.

%Also notice that in this case the sets $F^i_k(S_0)$ do not overlap, since they are contained in the $T_k^i$'s.
%Thus we can consider also the set $S:=\bigcup F_k^i(S_0)$. However,
%notice that this just amounts to consider the set $T'$ obtained by removing $T_0\setminus S_0$
%from the building block $T_0$.

%Now suppose that $\frac{\log b}{\log\lambda}\in(n-1,n)$. Then the set

\begin{rmk}
Roughly speaking, what matters is that there exists a bounded open set $T_0$ s.t.
\begin{equation*}
|F_k^i(T_0)\cap F_h^j(T_0)|=0,\qquad\textrm{if }i\not=j.
\end{equation*}
This can be thought of as a compatibility criterion for
the affine maps $\{F_k^i\}$.\\
We also need to ask that the ratio of the logarithms
of the growth factor and the scaling factor is $\frac{\log b}{\log\lambda}\in(n-1,n)$.\\
Then we are free to choose
as building block any set $T'_0\subset T_0$ s.t.
\begin{equation*}
|T'_0|>0,\qquad|T_0\setminus T'_0|>0\qquad\textrm{and }P(T'_0)<\infty,
\end{equation*}
and the set
\begin{equation*}
T':=\bigcup_{k=1}^\infty\bigcup_{i=1}^{ab^{k-1}}F_k^i(T'_0).
\end{equation*}
satisfies the hypothesis of previous Theorem.
\end{rmk}

Therefore, even if the Sierpinski triangle and the Menger sponge
do not satisfy $(\ref{add_frac_self_hp})$,
we can exploit their structure to construct new sets which do.

However, we remark that
the new boundary $\partial^-T'$ will look very different from the original fractal. Actually, in general
it will be a mix of unrectifiable pieces and smooth pieces. In particular, we can not hope to get an
analogue of Corollary $\ref{koch_coroll}$.
Still, the following Remark shows that the new (measure theoretic) boundary retains at least some of the ``fractal nature'' of the original set.

\begin{rmk}\label{self_sim_frac_bdry_nat_rmk}
If the set $T$ of Theorem $\ref{fractal_bdary_selfsim_dim}$ is bounded, exploiting Proposition $\ref{vis_prop}$ and Remark $\ref{fractal_dim_rmk}$
we obtain
\begin{equation*}
\overline{\Dim}_{\mathcal M}(\partial^-T)\geq\frac{\log b}{\log\lambda}>n-1.
\end{equation*}
Moreover, notice that if $\Omega$ is a bounded open set with Lipschitz boundary, then
\begin{equation}
P(E,\Omega)<\infty\quad\Longrightarrow\quad\Dim_F(E,\Omega)=n-1.
\end{equation}
Therefore, if $T\subset\subset B_R$, then
\begin{equation*}
P(T)=P(T,B_R)=\infty,
\end{equation*}
even if $T$ is bounded (and hence $\partial^-T$ is compact).

\end{rmk}

\begin{subsubsection}{Sponge-like sets}

The simplest way to construct the set $T'$ consists in simply removing a small ball
$S_0:=B\subset\subset T_0$ from $T_0$.

In particular, suppose that $|T_0\Delta T|=0$, as with the Sierpinski triangle.\\
Define
\begin{equation*}
S:=\bigcup_{k=1}^\infty\bigcup_{i=1}^{ab^{k-1}}F_k^i(B)
\quad\textrm{and}\quad
T':=\bigcup_{k=1}^\infty\bigcup_{i=1}^{ab^{k-1}}F_k^i(T_0\setminus B)=T\setminus S.
\end{equation*}
Then
\begin{equation}\label{fractal_spazz_end}
|T_0\Delta T|=0\quad\Longrightarrow\quad |T'\Delta (T_0\setminus S)|=0.
\end{equation}
Now the set $E:=T_0\setminus S$ looks like a sponge, in the sense that it
is a bounded open set with an infinite number of holes (each one at a positive, but non-fixed distance from the others).

From $(\ref{fractal_spazz_end})$ we get $P_s(E)=P_s(T')$. Thus, since $T'$ satisfies the hypothesis of
previous Theorem, we obtain
\begin{equation*}
\Dim_F(\partial^-E)=\frac{\log b}{\log\lambda}.
\end{equation*}

\end{subsubsection}

\begin{subsubsection}{Dendrite-like sets}

Depending on the form of the set $T_0$ and on the affine maps $\{F_k^i\}$,
we can define more intricated sets $T'$.

As an example we consider the Sierpinski triangle $E\subset\R^2$.\\
It is of the form $E=T_0\setminus T$,
where the building block $T_0$ is an equilateral triangle, say with side length one,
a vertex on the $y$-axis and baricenter in 0.
The pieces $T_k^i$ are obtained with a scaling factor $\lambda=2$
and the growth factor is $b=3$ (see e.g. \cite{Falconer} for the construction).
As usual, we consider the set
\begin{equation*}
T=\bigcup_{k=1}^\infty\bigcup_{i=1}^{3^{k-1}}T_k^i.
\end{equation*}
However, as remarked above, we have $|T\Delta T_0|=0$.

Starting from $k=2$ each triangle $T_k^i$ touches with (at least) a vertex (at least) another triangle $T_h^j$.
Moreover, each triangle $T_k^i$ gets touched in the middle point of each side (and actually it gets touched in infinitely many points).

Exploiting this situation, we can remove from $T_0$ six smaller triangles, so that the new building block
$T'_0$ is a star polygon centered in 0, with six vertices, one in each vertex of $T_0$ and one in each middle point of the sides of $T_0$.

\begin{figure}[htbp]
\begin{center}
\includegraphics[width=90mm]{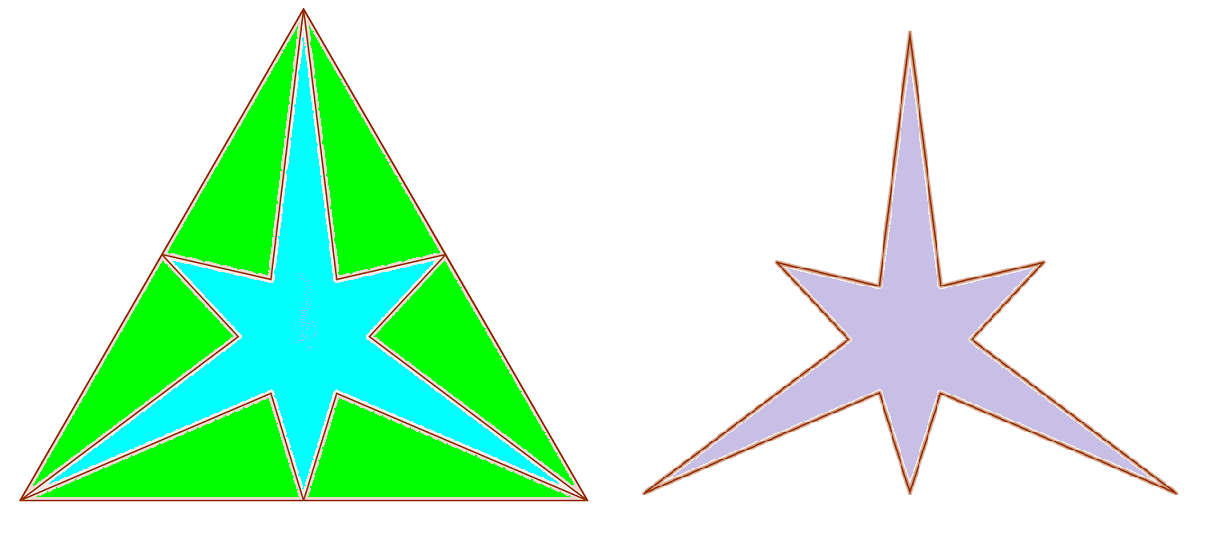}
\caption{{\it Removing the six triangles (in green) to obtain the new ``building block'' $T'_0$ (on the right)}}
\end{center}
\end{figure}

The resulting set
\begin{equation*}
T'=\bigcup_{k=1}^\infty\bigcup_{i=1}^{3^{k-1}}F_k^i(T'_0)
\end{equation*}
will have an infinite number of ramifications.

\begin{figure}[htbp]
\begin{center}
\includegraphics[width=110mm]{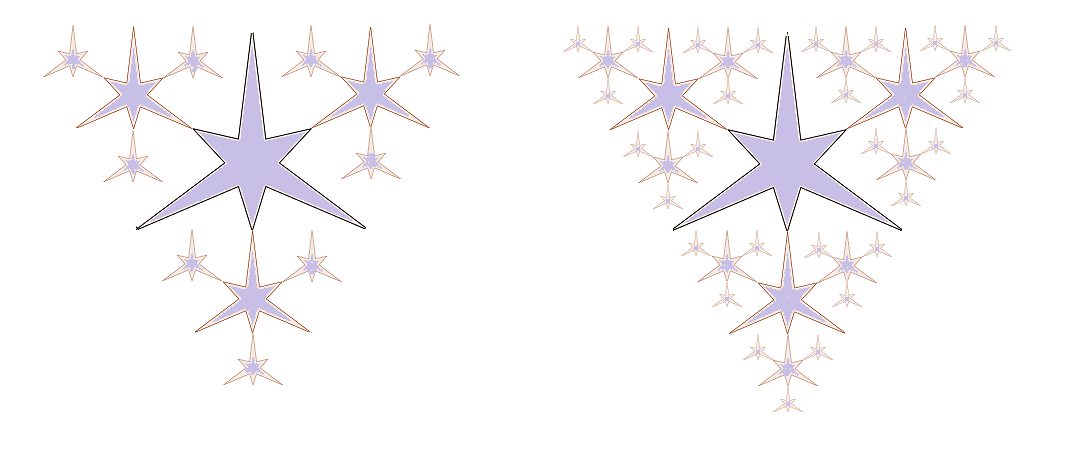}
\caption{{\it The third and fourth steps of the iterative construction of the set $T'$}}
\end{center}
\end{figure}

Since $T'$ satisfies the hypothesis of
previous Theorem, we obtain
\begin{equation*}
\Dim_F(\partial^-T')=\frac{\log 3}{\log2}.
\end{equation*}

\end{subsubsection}

\begin{subsubsection}{``Exploded'' fractals}

In all the previous examples, the sets $T_k^i$ are accumulated in a bounded region.

On the other hand, imagine making a fractal like the von Koch snowflake or the Sierpinski triangle
``explode'' and then rearrange the pieces $T_k^i$ in such a way that $d(T_k^i,T_h^j)\geq d$,
for some fixed $d>0$.

Since the shape of the building block is not important, we can consider $T_0:=B_{1/4}(0)\subset\R^n$,
with $n\geq2$. Moreover, since the parameter $a$ does not influence the dimension, we can fix $a=1$.

Then we rearrange the pieces obtaining
\begin{equation}\label{exploded_frac_def}
E:=\bigcup_{k=1}^\infty\bigcup_{i=1}^{b^{k-1}}B_\frac{1}{4\lambda^k}(k,0,\ldots,0,i).
\end{equation}
Define for simplicity
\begin{equation*}
B_k^i:=B_\frac{1}{4\lambda^k}(k,0,\ldots,0,i)\quad\textrm{and}\quad x_k^i:=k\,e_1+i\,e_n,
\end{equation*}
and notice that
\begin{equation*}
B_k^i=\lambda^{-k}B_\frac{1}{4}(0)+x_k^i.
\end{equation*}
Since for every $k,\,h$ and every $i\not=j$ we have
\begin{equation*}
d(B_k^i,B_h^j)\geq\frac{1}{2},
\end{equation*}
the boundary of the set $E$ is the disjoint union of $(n-1)$-dimensional spheres
\begin{equation*}
\partial^-E=\partial E=\bigcup_{k=1}^\infty\bigcup_{i=1}^{b^{k-1}}\partial B_k^i,
\end{equation*}
and in particular is smooth.

The (global) perimeter of $E$ is
\begin{equation*}
P(E)=\sum_{k=1}^\infty\sum_{i=1}^{b^{k-1}}P(B_k^i)=\frac{1}{\lambda}P(B_{1/4}(0))\sum_{k=0}^\infty
\Big(\frac{b}{\lambda^{n-1}}\Big)^k=\infty,
\end{equation*}
since $\frac{\log b}{\log\lambda}>n-1$.

However $E$ has locally finite perimeter, since its boundary is smooth and every ball $B_R$ intersects
only finitely many $B_k^i$'s,
\begin{equation*}
P(E,B_R)<\infty,\qquad\forall\,R>0.
\end{equation*}
Therefore it also has locally finite $s$-perimeter for every $s\in(0,1)$
\begin{equation*}
P_s(E,B_R)<\infty,\qquad\forall\,R>0,\qquad\forall\,s\in(0,1).
\end{equation*}

What is interesting is that the set $E$ satisfies the hypothesis of Theorem $\ref{fractal_bdary_selfsim_dim}$ and hence it also has finite global $s$-perimeter
for every $s<\sigma_0:=n-\frac{\log b}{\log\lambda}$,
\begin{equation*}
P_s(E)<\infty\qquad\forall\,s\in(0,\sigma_0)\quad\textrm{and}\quad P_s(E)=\infty\qquad\forall\,s\in[\sigma_0,1).
\end{equation*}

Thus we obtain Proposition $\ref{expl_farc_prop1}$.
%\begin{prop}\label{expl_farc_prop1}
%Let $n\geq2$. For every $\sigma\in(0,1)$ there exists a Caccioppoli set $E\subset\R^n$ s.t.
%\begin{equation}
%P_s(E)<\infty\qquad\forall\,s\in(0,\sigma)\quad\textrm{and}\quad P_s(E)=\infty\qquad\forall\,s\in[\sigma,1).
%\end{equation}

\begin{proof}[proof of Proposition $\ref{expl_farc_prop1}$]
It is enough to choose a natural number $b\geq2$ and take $\lambda:=b^\frac{1}{n-\sigma}$. Notice that
$\lambda>1$ and
\begin{equation*}
\frac{\log b}{\log\lambda}=n-\sigma\in(n-1,n).
\end{equation*}
Then we can define $E$ as in $(\ref{exploded_frac_def})$ and we are done.

\end{proof}
%\end{prop}

\end{subsubsection}

\end{subsection}

\begin{subsection}{Elementary properties of the $s$-perimeter}

\begin{prop}\label{elementary_properties}
Let $\Omega\subset\R^n$ be an open set.

(i) (Subadditivity)$\quad$ Let $E,\,F\subset\R^n$ s.t. $|E\cap F|=0$. Then
\begin{equation}\label{subadditive}
P_s(E\cup F,\Omega)\leq P_s(E,\Omega)+P_s(F,\Omega).
\end{equation}

(ii) (Translation invariance)$\quad$ Let $E\subset\R^n$ and $x\in\R^n$. Then
\begin{equation}\label{translation_invariance}
P_s(E+x,\Omega+x)=P_s(E,\Omega).
\end{equation}

(iii) (Rotation invariance)$\quad$ Let $E\subset\R^n$ and $\mathcal{R}\in SO(n)$ a rotation. Then
\begin{equation}\label{rotation_invariance}
P_s(\mathcal{R}E,\mathcal{R}\Omega)=P_s(E,\Omega).
\end{equation}

(iv) (Scaling)$\quad$ Let $E\subset\R^n$ and $\lambda>0$. Then
\begin{equation}\label{scaling}
P_s(\lambda E,\lambda\Omega)=\lambda^{n-s}P_s(E,\Omega).
\end{equation}
\begin{proof}
(i) follows from the following observations. Let $A_1,\,A_2,\,B\subset\R^n$. If $|A_1\cap A_2|=0$, then
\begin{equation*}
\Ll_s(A_1\cup A_2,B)
=\Ll_s(A_1,B)+\Ll_s(A_2,B).
\end{equation*}
Moreover
\begin{equation}
A_1\subset A_2\quad\Longrightarrow\quad\Ll_s(A_1,B)\leq\Ll_s(A_2,B),
\end{equation}
and
\begin{equation*}
\Ll_s(A,B)=\Ll_s(B,A).
\end{equation*}
Therefore
\begin{equation*}\begin{split}
P_s(E\cup F,\Omega)&=\Ll_s((E\cup F)\cap\Omega,\Co(E\cup F))+\Ll_s((E\cup F)\setminus\Omega,\Co(E\cup F)\cap\Omega)\\
&
=\Ll_s(E\cap\Omega,\Co(E\cup F))+\Ll_s(F\cap\Omega,\Co(E\cup F))\\
&
\qquad+\Ll_s(E\setminus\Omega,\Co(E\cup F)\cap\Omega)+\Ll_s(F\setminus\Omega,\Co(E\cup F)\cap\Omega)\\
&
\leq\Ll_s(E\cap\Omega,\Co E)+\Ll_s(F\cap\Omega,\Co F)\\
&
\qquad+\Ll_s(E\setminus\Omega,\Co E\cap\Omega)+\Ll_s(F\setminus\Omega,\Co F\cap\Omega)\\
&
=P_s(E,\Omega)+P_s(F,\Omega).
\end{split}\end{equation*}

(ii), (iii) and (iv) follow simply by changing variables in $\Ll_s$ and the following observations:
\begin{equation*}\begin{split}
&(x+A_1)\cap(x+A_2)=x+A_1\cap A_2,\qquad x+\Co A=\Co(x+A),\\
&
\mathcal{R}A_1\cap\mathcal{R}A_2=\mathcal{R}(A_1\cap A_2),\qquad\mathcal{R}(\Co A)=\Co(\mathcal{R}A),\\
&
(\lambda A_1)\cap(\lambda A_2)=\lambda(A_1\cap A_2),\qquad\lambda(\Co A)=\Co(\lambda A).
\end{split}\end{equation*}
For example, for claim (iv) we have
\begin{equation*}\begin{split}
\Ll_s(\lambda A,\lambda B)&=\int_{\lambda A}\int_{\lambda B}\frac{dx\,dy}{|x-y|^{n+s}}
=\int_A\lambda^n\,dx\int_B\frac{\lambda^n\,dy}{\lambda^{n+s}|x-y|^{n+s}}\\
&
=\lambda^{n-s}\Ll_s(A,B).
\end{split}
\end{equation*}
Then
\begin{equation*}\begin{split}
P_s(\lambda E,\lambda\Omega)&=\Ll_s(\lambda E\cap\lambda\Omega,\Co(\lambda E))+
\Ll_s(\lambda E\cap\Co(\lambda\Omega),\Co(\lambda E)\cap\lambda\Omega)\\
&
=\Ll_s(\lambda(E\cap\Omega),\lambda\Co E)+\Ll_s(\lambda(E\setminus\Omega),\lambda(\Co E\cap\Omega))\\
&
=\lambda^{n-s}\left(\Ll_s(E\cap\Omega,\Co E)+\Ll_s(E\setminus\Omega,\Co E\cap\Omega)\right)\\
&
=\lambda^{n-s}P_s(E,\Omega).
\end{split}\end{equation*}

\end{proof}
\end{prop}

\end{subsection}

\end{section}

\appendix

\section{Proof of Example $\ref{inclusion_counterexample}$}

Note that $E\subset (0,a^2]$.
Let $\Omega:=(-1,1)\subset\mathbb{R}$. Then $E\subset\subset\Omega$ and $\textrm{dist}(E,\partial\Omega)=1-a^2=:d>0$.
Now
\begin{equation*}
P_s(E)=\int_E\int_{\Co E\cap\Omega}\frac{dxdy}{|x-y|^{1+s}}+
\int_E\int_{\Co\Omega}\frac{dxdy}{|x-y|^{1+s}}.%=:I_1+I_2.
\end{equation*}
As for the second term, we have
\begin{equation*}
\int_E\int_{\Co\Omega}\frac{dxdy}{|x-y|^{1+s}}\leq\frac{2|E|}{sd^s}<\infty.
\end{equation*}
We split the first term into three pieces
\begin{equation*}\begin{split}
\int_E&\int_{\Co E\cap\Omega}\frac{dxdy}{|x-y|^{1+s}}\\
&
=\int_E\int_{-1}^0\frac{dxdy}{|x-y|^{1+s}}
+\int_E\int_{\Co E\cap(0,a)}\frac{dxdy}{|x-y|^{1+s}}+\int_E\int_a^1\frac{dxdy}{|x-y|^{1+s}}\\
&
=\mathcal{I}_1+\mathcal{I}_2+\mathcal{I}_3.
\end{split}
\end{equation*}
Note that $\Co E\cap(0,a)=\bigcup_{k\in\mathbb{N}}I_{2k-1}=\bigcup_{k\in\mathbb{N}}(a^{2k},a^{2k-1})$.\\
A simple calculation shows that, if $a<b\leq c<d$, then
\begin{equation}\label{rectangle_integral}\begin{split}
\int_a^b&\int_c^d\frac{dxdy}{|x-y|^{1+s}}=\\
&
\frac{1}{s(1-s)}\big[(c-a)^{1-s}+(d-b)^{1-s}-(c-b)^{1-s}-(d-a)^{1-s}\big].
\end{split}
\end{equation}
Also note that, if $n>m\geq1$, then
\begin{equation}\label{derivative_bound}\begin{split}
(1-a^n)^{1-s}-(1-a^m)^{1-s}&=\int_m^n\frac{d}{dt}(1-a^t)^{1-s}\,dt\\
&
=(s-1)\log a\int_m^n\frac{a^t}{(1-a^t)^s}\,dt\\
&
\leq a^m (s-1)\log a\int_m^n\frac{1}{(1-a^t)^s}\,dt\\
&
\leq(n-m)a^m\frac{(s-1)\log a}{(1-a)^s}.
\end{split}
\end{equation}
Now consider the first term
\begin{equation*}
\mathcal{I}_1=\sum_{k=1}^\infty\int_{a^{2k+1}}^{a^{2k}}\int_{-1}^0\frac{dxdy}{|x-y|^{1+s}}.
\end{equation*}
Use $(\ref{rectangle_integral}$) and notice that $(c-a)^{1-s}-(d-a)^{1-s}\leq0$ to get
\begin{equation*}
\int_{-1}^0\int_{a^{2k+1}}^{a^{2k}}\frac{dxdy}{|x-y|^{1+s}}
\leq\frac{1}{s(1-s)}\big[(a^{2k})^{1-s}-(a^{2k+1})^{1-s}\big]\leq\frac{1}{s(1-s)}(a^{2(1-s)})^k.
\end{equation*}
Then, as $a^{2(1-s)}<1$ we get
\begin{equation*}
\mathcal{I}_1\leq\frac{1}{s(1-s)}\sum_{k=1}^\infty(a^{2(1-s)})^k<\infty.
\end{equation*}
As for the last term
\begin{equation*}
\mathcal{I}_3=\sum_{k=1}^\infty\int_{a^{2k+1}}^{a^{2k}}\int_a^1\frac{dxdy}{|x-y|^{1+s}},
\end{equation*}
use $(\ref{rectangle_integral}$) and notice that $(d-b)^{1-s}-(d-a)^{1-s}\leq0$ to get
\begin{equation*}\begin{split}
\int_{a^{2k+1}}^{a^{2k}}\int_a^1\frac{dxdy}{|x-y|^{1+s}}&
\leq\frac{1}{s(1-s)}\big[(1-a^{2k+1})^{1-s}-(1-a^{2k})^{1-s}\big]\\
&
\leq\frac{-\log a}{s(1-a)^s}a^{2k}\quad\textrm{by }(\ref{derivative_bound}).
\end{split}
\end{equation*}
Thus
\begin{equation*}
\mathcal{I}_3\leq\frac{-\log a}{s(1-a)^s}\sum_{k=1}^\infty(a^2)^k<\infty.
\end{equation*}
Finally we split the second term
\begin{equation*}
\mathcal{I}_2=\sum_{k=1}^\infty\sum_{j=1}^\infty\int_{a^{2k+1}}^{a^{2k}}\int_{a^{2j}}^{a^{2j-1}}
\frac{dxdy}{|x-y|^{1+s}}
\end{equation*}
into three pieces according to the cases $j>k$, $j=k$ and $j<k$.

If $j=k$, using $(\ref{rectangle_integral})$ we get
\begin{equation*}\begin{split}
\int_{a^{2k+1}}^{a^{2k}}&\int_{a^{2k}}^{a^{2k-1}}
\frac{dxdy}{|x-y|^{1+s}}=\\
&
=\frac{1}{s(1-s)}\big[(a^{2k}-a^{2k+1})^{1-s}+(a^{2k-1}-a^{2k})^{1-s}-(a^{2k-1}-a^{2k+1})^{1-s}\big]\\
&
=\frac{1}{s(1-s)}\big[a^{2k(1-s)}(1-a)^{1-s}+a^{(2k-1)(1-s)}(1-a)^{1-s}\\
&
\quad\quad\quad\quad\quad-a^{(2k-1)(1-s)}(1-a^2)^{1-s}\big]\\
&
=\frac{1}{s(1-s)}(a^{2(1-s)})^k\Big[(1-a)^{1-s}+\frac{(1-a)^{1-s}}{a^{1-s}}-\frac{(1-a^2)^{1-s}}{a^{1-s}}\Big].
\end{split}
\end{equation*}
Summing over $k\in\mathbb{N}$ we get
\begin{equation*}\begin{split}
\sum_{k=1}^\infty&\int_{a^{2k+1}}^{a^{2k}}\int_{a^{2k}}^{a^{2k-1}}
\frac{dxdy}{|x-y|^{1+s}}=\\
&
=\frac{1}{s(1-s)}\frac{a^{2(1-s)}}{1-a^{2(1-s)}}\Big[(1-a)^{1-s}+\frac{(1-a)^{1-s}}{a^{1-s}}-\frac{(1-a^2)^{1-s}}{a^{1-s}}\Big]<\infty.
\end{split}
\end{equation*}
In particular note that
\begin{equation*}\begin{split}
(1-s)&P_s(E)\geq(1-s)\mathcal{I}_2\\
&
\geq\frac{1}{s(1-a^{2(1-s)})}\big[a^{2(1-s)}(1-a)^{1-s}+a^{1-s}(1-a)^{1-s}-a^{1-s}(1-a^2)^{1-s}\big],
\end{split}
\end{equation*}
which tends to $+\infty$ when $s\to1$. This shows that $E$ cannot have finite perimeter.

To conclude let $j>k$, the case $j<k$ being similar, and consider
\begin{equation*}
\sum_{k=1}^\infty\sum_{j=k+1}^\infty\int_{a^{2j}}^{a^{2j-1}}\int_{a^{2k+1}}^{a^{2k}}
\frac{dxdy}{|x-y|^{1+s}}.
\end{equation*}
Again, using $(\ref{rectangle_integral}$) and $(d-b)^{1-s}-(d-a)^{1-s}\leq0$, we get
\begin{equation*}\begin{split}
\int_{a^{2j}}^{a^{2j-1}}&\int_{a^{2k+1}}^{a^{2k}}
\frac{dxdy}{|x-y|^{1+s}}\\
&
\leq\frac{1}{s(1-s)}\big[(a^{2k+1}-a^{2j})^{1-s}-(a^{2k+1}-a^{2j-1})^{1-s}\big]\\
&
=\frac{a^{1-s}}{s(1-s)}(a^{2(1-s)})^k\big[(1-a^{2(j-k)-1})^{1-s}-(1-a^{2(j-k)-2})^{1-s}\big]\\
&
\leq\frac{a^{1-s}}{s(1-s)}(a^{2(1-s)})^k\frac{(s-1)\log a}{(1-a)^s}a^{2(j-k)-2}\quad\quad\textrm{by }(\ref{derivative_bound})\\
&
=\frac{-\log a}{s(1-a^s)a^{s+1}}(a^{2(1-s)})^k(a^2)^{j-k},
\end{split}
\end{equation*}
for $j\geq k+2$. Then
\begin{equation*}
\begin{split}
\sum_{k=1}^\infty&\sum_{j=k+2}^\infty\int_{a^{2j}}^{a^{2j-1}}\int_{a^{2k+1}}^{a^{2k}}
\frac{dxdy}{|x-y|^{1+s}}\\
&
\leq\frac{-\log a}{s(1-a^s)a^{s+1}}\sum_{k=1}^\infty(a^{2(1-s)})^k\sum_{h=2}^\infty(a^2)^h<\infty.
\end{split}
\end{equation*}
If $j=k+1$ we get
\begin{equation*}\begin{split}
\sum_{k=1}^\infty\int_{a^{2k+2}}^{a^{2k+1}}\int_{a^{2k+1}}^{a^{2k}}\frac{dxdy}{|x-y|^{1+s}}&
\leq\frac{1}{s(1-s)}\sum_{k=1}^\infty(a^{2k+1}-a^{2k+2})^{1-s}\\
&
=\frac{a^{1-s}(1-a)^{1-s}}{s(1-s)}\sum_{k=1}^\infty(a^{2(1-s)})^k<\infty.
\end{split}
\end{equation*}
This shows that also $\mathcal{I}_2<\infty$, so that $P_s(E)<\infty$ for every $s\in(0,1)$ as claimed.

\section{Signed distance function}

Given $\emptyset\not=E\subset\R^n$, the distance function from $E$ is defined as
\begin{equation*}
d_E(x)=d(x,E):=\inf_{y\in E}|x-y|,\qquad\textrm{for }x\in\R^n.
\end{equation*}
The signed distance function from $\partial E$, negative inside $E$, is then defined as
\begin{equation}
\bar{d}_E(x)=\bar{d}(x,E):=d(x,E)-d(x,\Co E).
\end{equation}
For the details of the main properties we refer e.g. to \cite{Ambrosio} and \cite{Bellettini}.

We also define the sets
\begin{equation*}
E_r:=\{x\in\R^n\,|\,\bar{d}_E(x)<r\}.
\end{equation*}

Let $\Omega\subset\R^n$ be a bounded open set with Lipschitz boundary. By definition we can locally describe $\Omega$ near its boundary as the subgraph of appropriate Lipschitz functions.
To be more precise, we can find a finite open covering $\{C_{\rho_i}\}_{i=1}^m$ of $\partial\Omega$ made of cylinders,
and Lipschitz functions $\varphi_i:B'_{\rho_i}\longrightarrow\R$ s.t. $\Omega\cap C_{\rho_i}$ is the subgraph of
$\varphi_i$.
That is, up to rotations and translations,
\begin{equation*}
C_{\rho_i}=\{(x',x_n)\in\R^n\,|\,|x'|<\rho_i,\,|x_n|<\rho_i\},
\end{equation*}
and
\begin{equation*}\begin{split}
\Omega\cap C_{\rho_i}&=\{(x',x_n)\in\R^n\,|\,x'\in B'_{\rho_i},\,-\rho_i<x_n<\varphi_i(x')\},\\
&
\partial\Omega\cap C_{\rho_i}=\{(x',\varphi_i(x'))\in\R^n\,|\,x'\in B_{\rho_i}'\}.
\end{split}
\end{equation*}
Let $L$ be the sup of the Lipschitz constants of the functions $\varphi_i$.

Theorem 4.1 of \cite{LipApprox} guarantees that also the bounded open sets $\Omega_r$ have
Lipschitz boundary, when $r$ is small enough, say $|r|<r_0$.\\
Moreover these sets $\Omega_r$
can locally be described, in the same cylinders $C_{\rho_i}$ used for $\Omega$, as subgraphs of
Lipschitz functions $\varphi_i^r$ which approximate $\varphi_i$ (see \cite{LipApprox} for the precise statement) and whose Lipschitz constants are less or equal to $L$.\\
Notice that
\begin{equation*}
\partial\Omega_r=\{\bar{d}_\Omega=r\}.
\end{equation*}
Now, since in $C_{\rho_i}$ the set $\Omega_r$ coincides with the subgraph of $\varphi_i^r$, we have
\begin{equation*}
\Ha^{n-1}(\partial\Omega_r\cap C_{\rho_i})=\int_{B_{\rho_i}'}\sqrt{1+|\nabla\varphi_i^r|^2}\,dx'\leq M_i,
\end{equation*}
with $M_i$ depending on $\rho_i$ and $L$ but not on $r$.\\
Therefore
\begin{equation*}
\Ha^{n-1}(\{\bar{d}_\Omega=r\})\leq\sum_{i=1}^m\Ha^{n-1}(\partial\Omega_r\cap C_{\rho_i})\leq\sum_{i=1}^mM_i
\end{equation*}
independently on $r$,
proving the following
\begin{prop}\label{bound_perimeter_unif}
Let $\Omega\subset\R^n$ be a bounded open set with Lipschitz boundary. Then there exists $r_0>0$ s.t.
$\Omega_r$ is a bounded open set with Lipschitz boundary for every $r\in(-r_0,r_0)$ and
\begin{equation}\label{bound_perimeter_unif_eq}
\sup_{|r|<r_0}\Ha^{n-1}(\{\bar{d}_\Omega=r\})<\infty.
\end{equation}
\end{prop}

\section{Measure theoretic boundary}

Since
\begin{equation}\label{fin_spazz_basta1}
|E\Delta F|=0\quad\Longrightarrow\quad P(E,\Omega)=P(F,\Omega)\quad\textrm{and}\quad P_s(E,\Omega)=P_s(F,\Omega),
\end{equation}
we can modify a set making its topological boundary as big as we want, without changing its (fractional) perimeter.\\
For example, let $E\subset\R^n$ be a bounded open set with Lipschitz boundary. Then, if we set
$F:=E\cup(\mathbb Q^n\setminus E)$, we have $|E\Delta F|=0$ and hence we get $(\ref{fin_spazz_basta1})$.
However $\partial F=\R^n\setminus E$.

For this reason one considers measure theoretic notions of interior, exterior and boundary, which solely depend on the class
of $\chi_E$ in $L^1_{loc}(\R^n)$.\\
In some sense, by considering the measure theoretic boundary $\partial^-E$ defined below
we can also minimize the size of the topological boundary (see $(\ref{ess_bdry_intersect})$). Moreover, this measure theoretic boundary is actually the
topological boundary of a set which is equivalent to $E$. Thus we obtain a ``good'' representative for the class of $E$.

We refer to Section 3.2 of \cite{Visintin} (see also Proposition 3.1 of \cite{Giusti}). For some details about the good representative of an $s$-minimal set, see the Appendix of \cite{graph}.
\begin{defin}
Let $E\subset\R^n$. For every $t\in[0,1]$ define the set
\begin{equation}\label{density_t}
E^{(t)}:=\left\{x\in\R^n\,\big|\,\exists\lim_{r\to0}\frac{|E\cap B_r(x)|}{\omega_nr^n}=t\right\},
\end{equation}
of points density $t$ of $E$. The sets $E^{(0)}$ and $E^{(1)}$ are respectively the measure theoretic exterior and interior of the set $E$. The set
\begin{equation}\label{ess_bdry}
\partial_eE:=\R^n\setminus(E^{(0)}\cup E^{(1)})
\end{equation}
is the essential boundary of $E$.
\end{defin}

Using the Lebesgue points Theorem for the characteristic function $\chi_E$, we see that the limit in $(\ref{density_t})$ exists for a.e. $x\in\R^n$ and
\begin{equation*}
\lim_{r\to0}\frac{|E\cap B_r(x)|}{\omega_nr^n}=\left\{\begin{array}{cc}1,&\textrm{a.e. }x\in E,\\
0,&\textrm{a.e. }x\in\Co E.
\end{array}
\right.
\end{equation*}
So
\begin{equation*}
|E\Delta E^{(1)}|=0,\qquad|\Co E\Delta E^{(0)}|=0\qquad\textrm{and }|\partial_eE|=0.
\end{equation*}
In particular every set $E$ is equivalent to its measure theoretic interior.\\
However, notice that $E^{(1)}$ in general is not open.\\

We have another natural way to define a measure theoretic boundary.
\begin{defin}
Let $E\subset\R^n$ and define the sets
\begin{equation*}\begin{split}
&E_1:=\{x\in\R^n\,|\,\exists r>0,\,|E\cap B_r(x)|=\omega_nr^n\},\\
&
E_0:=\{x\in\R^n\,|\,\exists r>0,\,|E\cap B_r(x)|=0\}.
\end{split}\end{equation*}
Then we define
\begin{equation*}\begin{split}
\partial^-E&:=\R^n\setminus(E_0\cup E_1)\\
&
=\{x\in\R^n\,|\,0<|E\cap B_r(x)|<\omega_nr^n\textrm{ for every }r>0\}.
\end{split}
\end{equation*}
\end{defin}
Notice that $E_0$ and $E_1$ are open sets and hence $\partial^-E$ is closed. Moreover, since
\begin{equation}\label{density_subsets}
E_0\subset E^{(0)}\qquad\textrm{and}\qquad E_1\subset E^{(1)},
\end{equation}
we get
\begin{equation*}
\partial_eE\subset\partial^-E.
\end{equation*}
We have
\begin{equation}\label{ess_bdry_top1}
F\subset\R^n\textrm{ s.t. }|E\Delta F|=0\quad\Longrightarrow\quad\partial^-E\subset\partial F.
\end{equation}
Indeed, if $|E\Delta F|=0$, then $|F\cap B_r(x)|=|E\cap B_r(x)|$ for every $r>0$. Thus for any $x\in\partial^-E$ we have
\begin{equation*}
0<|F\cap B_r(x)|<\omega_nr^n,
\end{equation*}
which implies
\begin{equation*}
F\cap B_r(x)\not=\emptyset\quad\textrm{and}\quad\Co F\cap B_r(x)\not=\emptyset\quad\textrm{for every }r>0,
\end{equation*}
and hence $x\in\partial F$.

In particular, $\partial^-E\subset\partial E$.

Moreover
\begin{equation}\label{ess_bdry_top2}
\partial^-E=\partial E^{(1)}.
\end{equation}
Indeed, since $|E\Delta E^{(1)}|=0$, we already know that $\partial^-E\subset\partial E^{(1)}$.
The converse inclusion follows from $(\ref{density_subsets})$ and the fact that both $E_0$ and $E_1$ are open.\\
From $(\ref{ess_bdry_top1})$ and $(\ref{ess_bdry_top2})$ we obtain
\begin{equation}\label{ess_bdry_intersect}
\partial^-E=\bigcap_{F\sim E}\partial F,
\end{equation}
where the intersection is taken over all sets $F\subset\R^n$ s.t. $|E\Delta F|=0$,
so we can think of 
$\partial^-E$ as a way to minimize the size of the topological boundary of $E$.
In particular
\begin{equation*}
F\subset\R^n\textrm{ s.t. }|E\Delta F|=0\quad\Longrightarrow\quad\partial^-F=\partial^-E.
\end{equation*}

From $(\ref{density_subsets})$ and $(\ref{ess_bdry_top2})$ we see that we can take $E^{(1)}$ as ``good'' representative
for $E$, obtaining Remark $\ref{gmt_assumption}$.\\

Recall that the support of a Radon measure $\mu$ on $\R^n$ is defined as the set
\begin{equation*}
\textrm{supp }\mu:=\{x\in\R^n\,|\,\mu(B_r(x))>0\textrm{ for every }r>0\}.
\end{equation*}
Notice that, being the complementary of the union of all open sets of measure zero, it is a closed set.
In particular, if $E$ is a Caccioppoli set, we have
\begin{equation}\label{support_perimeter}
\textrm{supp }|D\chi_E|=\{x\in\R^n\,|\,P(E,B_r(x))>0\textrm{ for every }r>0\},
\end{equation}
and
it is easy to verify that
\begin{equation*}
\partial^-E=\textrm{supp }|D\chi_E|=\overline{\partial^*E},
\end{equation*}
where $\partial^*E$ denotes the reduced boundary.
However notice that in general the inclusions
\begin{equation*}
\partial^*E\subset\partial_eE\subset\partial^-E\subset\partial E
\end{equation*}
are all strict and in principle we could have
\begin{equation*}
\Ha^{n-1}(\partial^-E\setminus\partial^*E)>0.
\end{equation*}

\section{Minkowski dimension}

\begin{defin}
Let $\Omega\subset\R^n$ be an open set. For any $\Gamma\subset\R^n$ and $r\in[0,n]$ we define the inferior and superior $r$-dimensional Minkowski contents of $\Gamma$ relative to the set $\Omega$ as, respectively
\begin{equation*}
\underline{\mathcal{M}}^r(\Gamma,\Omega):=\liminf_{\rho\to0}\frac{|\bar{N}_\rho^\Omega(\Gamma)|}{\rho^{n-r}},\qquad
\overline{\mathcal{M}}^r(\Gamma,\Omega):=\limsup_{\rho\to0}\frac{|\bar{N}_\rho^\Omega(\Gamma)|}{\rho^{n-r}}.
\end{equation*}
Then we define the lower and upper Minkowski dimensions of $\Gamma$ in $\Omega$ as
\begin{equation*}\begin{split}
\underline{\Dim}_\mathcal{M}(\Gamma,\Omega)&:=\inf\big\{r\in[0,n]\,|\,\underline{\mathcal{M}}^r(\Gamma,\Omega)=0\big\}\\
&
=n-\sup\big\{r\in[0,n]\,|\,\underline{\mathcal{M}}^{n-r}(\Gamma,\Omega)=0\big\},
\end{split}\end{equation*}
\begin{equation*}\begin{split}
\overline{\Dim}_\mathcal{M}(\Gamma,\Omega)&:=\sup\big\{r\in[0,n]\,|\,\overline{\mathcal{M}}^r(\Gamma,\Omega)=\infty\big\}\\
&
=n-\inf\big\{r\in[0,n]\,|\,\overline{\mathcal{M}}^{n-r}(\Gamma,\Omega)=\infty\big\}.
\end{split}
\end{equation*}
If they agree, we write
\begin{equation*}
\Dim_\mathcal{M}(\Gamma,\Omega)
\end{equation*}
for the common value and call it the Minkowski dimension of $\Gamma$ in $\Omega$.\\
If $\Omega=\R^n$ or $\Gamma\subset\subset\Omega$, we drop the $\Omega$ in the formulas.
\end{defin}

\begin{rmk}
Let $\Dim_\mathcal{H}$ denote the Hausdorff dimension. In general one has
\begin{equation*}
\Dim_\mathcal{H}(\Gamma)\leq\underline{\Dim}_\mathcal{M}(\Gamma)\leq\overline{\Dim}_\mathcal{M}(\Gamma),
\end{equation*}
and all the inequalities might be strict. However for some sets (e.g. self-similar sets with some symmetric and regularity condition) they are all equal.
\end{rmk}

\noindent
We also recall some equivalent definitions of the Minkowski dimensions, usually referred to as box-counting dimensions, which are easier to compute. For the details and the relation between the Minkowski and the Hausdorff dimensions, see \cite{Mattila} and \cite{Falconer} and the references cited therein.\\
For simplicity we only consider the case $\Gamma$ bounded and $\Omega=\R^n$ (or $\Gamma\subset\subset\Omega$).
\begin{defin}
Given a nonempty bounded set $\Gamma\subset\R^n$, define for every $\delta>0$
\begin{equation*}
\mathcal{N}(\Gamma,\delta):=\min\Big\{k\in\mathbb{N}\,\big|\,\Gamma\subset\bigcup_{i=1}^kB_\delta(x_i),\textrm{ for some }x_i\in\R^n\Big\},
\end{equation*}
the smallest number of $\delta$-balls needed to cover $\Gamma$, and
\begin{equation*}
\mathcal{P}(\Gamma,\delta):=\max\big\{k\in\mathbb{N}\,|\,\exists\,\textrm{disjoint balls }B_\delta(x_i),\,i=1,\ldots,k\textrm{ with }x_i\in \Gamma\big\},
\end{equation*}
the greatest number of disjoint $\delta$-balls with centres in $\Gamma$.
\end{defin}

Then it is easy to verify that
\begin{equation}\label{counting}
\mathcal{N}(\Gamma,2\delta)\leq\mathcal{P}(\Gamma,\delta)\leq\mathcal{N}(\Gamma,\delta/2).
\end{equation}
Moreover, since any union of $\delta$-balls with centers in $\Gamma$ is contained in $N_\delta(\Gamma)$, and any
union of $(2\delta)$-balls covers $N_\delta(\Gamma)$ if the union of the corresponding $\delta$-balls covers $\Gamma$, we get
\begin{equation}\label{counting2}
\mathcal{P}(\Gamma,\delta)\omega_n\delta^n\leq|N_\delta(\Gamma)|\leq
\mathcal{N}(\Gamma,\delta)\omega_n(2\delta)^n.
\end{equation}
Using $(\ref{counting})$ and $(\ref{counting2})$ we see that
\begin{equation*}\begin{split}
&\underline{\Dim}_\mathcal{M}(\Gamma)=\inf\Big\{r\in[0,n]\,\big|\,\liminf_{\delta\to0}\mathcal{N}(\Gamma,\delta)\delta^r=0\Big\},\\
&
\overline{\Dim}_\mathcal{M}(\Gamma)=\sup\Big\{r\in[0,n]\,\big|\,\limsup_{\delta\to0}\mathcal{N}(\Gamma,\delta)\delta^r=\infty\Big\}.
\end{split}
\end{equation*}
Then it can be proved that
\begin{equation}\label{log_counting}\begin{split}
&\underline{\Dim}_\mathcal{M}(\Gamma)=\liminf_{\delta\to0}\frac{\log\mathcal{N}(\Gamma,\delta)}{-\log\delta},\\
&
\overline{\Dim}_\mathcal{M}(\Gamma)=\limsup_{\delta\to0}\frac{\log\mathcal{N}(\Gamma,\delta)}{-\log\delta}.
\end{split}
\end{equation}
Actually notice that, due to $(\ref{counting})$,
we can take $\mathcal{P}(\Gamma,\delta)$ in place of $\mathcal{N}(\Gamma,\delta)$ in the above formulas.\\
It is also easy to see that if in the definition of $\mathcal{N}(\Gamma,\delta)$ we take cubes of side $\delta$ instead of balls of radius $\delta$, then we get exactly the same dimensions.

Moreover in $(\ref{log_counting})$ it is enough to consider limits as $\delta\to0$
through any decreasing sequence $\delta_k$
s.t. $\delta_{k+1}\geq c\delta_k$ for some constant $c\in(0,1)$; in particular for $\delta_k=c^k$. Indeed
if $\delta_{k+1}\leq\delta<\delta_k$, then
\begin{equation*}\begin{split}
\frac{\log\mathcal{N}(\Gamma,\delta)}{-\log\delta}&\leq\frac{\log\mathcal{N}(\Gamma,\delta_{k+1})}{-\log\delta_k}
=\frac{\log\mathcal{N}(\Gamma,\delta_{k+1})}{-\log\delta_{k+1}+\log(\delta_{k+1}/\delta_k)}\\
&
\leq\frac{\log\mathcal{N}(\Gamma,\delta_{k+1})}{-\log\delta_{k+1}+\log c},
\end{split}\end{equation*}
so that
\begin{equation*}
\limsup_{\delta\to0}\frac{\log\mathcal{N}(\Gamma,\delta)}{-\log\delta}\leq
\limsup_{k\to\infty}\frac{\log\mathcal{N}(\Gamma,\delta_k)}{-\log\delta_k}.
\end{equation*}
The opposite inequality is clear and in a similar way we can treat the lower limits.

\end{document}